\documentclass[10pt]{amsart}

\usepackage{graphicx}
\usepackage{amsthm}
\usepackage{datetime}
\usepackage{mathbbol}

\newcommand{\sujap}{{\mathfrak S}}
\newcommand{\sajap}{{\mathfrak M}}

\newcommand{\noi}{\noindent}

\newcommand{\ltrbr}{[\mkern-3mu[\mkern-3mu[}

\newcommand{\rtrbr}{]\mkern-3mu]\mkern-3mu]}

\usepackage{textcomp}
\usepackage{mathabx}
\input bookman.sty
\boldmath

\usepackage{datetime}

\usepackage{comment}
\excludecomment{excl}

\usepackage{xcolor}

\usepackage{helvet}

\usepackage{graphics}
\usepackage{amssymb}
\usepackage{amsxtra}
\usepackage{amsmath}
\usepackage{mathrsfs}

\usepackage[arrow, matrix]{xy}
\xyoption{frame}

\theoremstyle{plain}

\newtheorem{theo}{Theorem}[section]
\newtheorem{lemm}[theo]{Lemma}
\newtheorem{prop}[theo]{Proposition}
\newtheorem{coro}[theo]{Corollary}

\theoremstyle{definition}

\newtheorem{defi}[theo]{Definition}
\newtheorem{rema}[theo]{Remark}

\newtheorem{exam}[theo]{Example}

\newfont{\rmm}{cmr10 scaled 1000}

\newfont{\itt}{cmsl10 scaled 1000}

\newfont{\rM}{cmr10 scaled 1700}

\newcounter{lemma}[section]

\newcounter{tempcounter}

\newcommand{\lb}{\label}

\newcommand{\rrf}[1]{(\ref{#1})}

\renewcommand{\a}{\alpha}
\renewcommand{\b}{\beta}
\newcommand{\g}{\gamma}
\renewcommand{\d}{\delta}
\newcommand{\e}{\epsilon}
\newcommand{\ve}{\varepsilon}
\newcommand{\z}{\zeta}
\renewcommand{\t}{\theta}
\renewcommand{\l}{\lambda}
\renewcommand{\k}{\varkappa}
\newcommand{\m}{\mu}

\renewcommand{\r}{\rho}

\newcommand{\s}{\sigma}

\renewcommand{\o}{\omega}

\newcommand{\G}{\Gamma}
\newcommand{\D}{\Delta}

\renewcommand{\L}{\Lambda}

\renewcommand{\O}{\Omega}

\newcommand{\BB}{{\mathcal B}}
\newcommand{\CC}{{\mathcal C}}

\newcommand{\FF}{{\mathcal F}}
\newcommand{\GG}{{\mathcal G}}
\newcommand{\HH}{{\mathcal H}}

\newcommand{\NN}{{\mathcal N}}

\renewcommand{\SS}{{\mathcal S}}

\newcommand{\UU}{{\mathcal U}}

\newcommand{\nn}{{\mathbb{N}}}

\newcommand{\qq}{{\mathbb{Q}}}
\newcommand{\rr}{{\mathbb{R}}}

\newcommand{\ttt}{{\mathbb{T}}}

\newcommand{\zz}{{\mathbb{Z}}}

\newcommand{\AAAA}{{\mathscr{A}}}
\newcommand{\BBBB}{{\mathscr{B}} }

\newcommand{\DDDD}{{\mathscr{D}} }
\newcommand{\EEEE}{{\mathscr{E}} }
\newcommand{\FFFF}{{\mathscr{F}} }

\newcommand{\ZZZZ}{{\mathscr{Z}}}

\newcommand{\gC}{{\mathfrak{C}}}

\newcommand{\id}{\text{id}}

\newcommand{\Ker}{\text{\rm Ker }}

\newcommand{\Hom}{\text{\rm Hom}}

\newcommand{\ind}{\text{\rm ind\hspace{0.05cm}}}

\renewcommand{\Im}{\text{\rm Im }}
\newcommand{\supp}{\text{\rm supp }}
\newcommand{\Int}{\text{\rm Int }}

\newcommand{\Id}{\text{\rm Id}}

\newcommand{\bere}{\begin{rema}}
\newcommand{\bede}{\begin{defi}}

\renewcommand{\beth}{\begin{theo}}
\newcommand{\bele}{\begin{lemm}}
\newcommand{\bepr}{\begin{prop}}
\newcommand{\beeq}{\begin{equation}}
\newcommand{\bega}{\begin{gather}}
\newcommand{\begaa}{\begin{gather*}}
\newcommand{\been}{\begin{enumerate}}

\newcommand{\bedee}{\begin{defii}}
\newcommand{\bethh}{\begin{theoo}}
\newcommand{\belee}{\begin{lemmm}}
\newcommand{\beprr}{\begin{propp}}

\newcommand{\beco}{\begin{coro}}

\newcommand{\beal}{\begin{aligned}}

\newcommand{\enre}{\end{rema}}

\newcommand{\enco}{\end{coro}}
\newcommand{\enpr}{\end{prop}}
\newcommand{\enth}{\end{theo}}
\newcommand{\enle}{\end{lemm}}
\newcommand{\enen}{\end{enumerate}}
\newcommand{\enga}{\end{gather}}
\newcommand{\engaa}{\end{gather*}}
\newcommand{\eneq}{\end{equation}}
\newcommand{\enal}{\end{aligned}}

\newcommand{\bq}{\begin{equation}}
\newcommand{\bqq}{\begin{equation*}}

\renewcommand{\leq}{\leqslant}
\renewcommand{\geq}{\geqslant}

\newcommand{\mxx}[1]{\quad\mbox{#1}\quad}
 
\newcommand{\wi}{\widetilde}

\newcommand{\ove}{\overline}

\newcommand{\wh}{\widehat}

\newcommand{\sm}{\setminus}

\newcommand{\sbs}{\subset}

\newcommand{\str}[1]{{#1}^{\displaystyle\twoheadrightarrow}}

\newcommand{\st}[1]{\overset{\rightsquigarrow}{#1}}

\newcommand{\stexp}[1]{{#1}^{\rightsquigarrow}}

\newcommand{\stv}{\stexp {(-v)}}

\newcommand{\tens}[1]{\underset{#1}{\otimes}}

\newcommand{\Lxi}{{\wh \L}_\xi}

\newcommand{\lL}{\wh{\wh L}}

\newcommand{\dow}{\pr_0 W}

\newcommand{\daw}{\pr_1 W}

\newcommand{\sut}{~such~that~}

\newcommand{\wrt}{with respect to}
\newcommand{\ho}{homomorphism}
\newcommand{\iso}{isomorphism}

\newcommand{\ma}{manifold}
\newcommand{\nei}{neighbourhood}

\newcommand{ \co}{~cobordism}
\newcommand{
\sma}{submanifold}

\newcommand{\noconf}{~no~confusion~is~possible}

\newcommand{\mfcob}{~ Let $\fcob$
be a Morse function on
 a cobordism $W$}

\newcommand{\Prf}{{\it Proof.\quad}}
\newcommand{\prf}{{\it Proof:\quad}}

\newcommand{\smo}{C^{\infty}}

\newcommand{\chart}{\Phi_p:U_p\to B^n(0,r_p)}
\newcommand{\atlas}{\{\Phi_p:U_p\to B^n(0,r_p)\}_{p\in S(f)}}

\newcommand{\fcob}{f:W\to[a,b]}

\newcommand{\indl}[1]{{\scriptstyle{\text{\rm ind}\leqslant {#1}~}}}

\newcommand{\pr}{\partial}

\newcommand{\qt}{\hfill\triangle}
\newcommand{\qs}{\hfill\square}

\newcommand{\pa}{\vskip0.1in}

\newcommand{\liminv}{\underset {\leftarrow}{\lim}}

\newcommand{\Lxim}{{\wh \L}^-_\xi}

\newcommand{\lLL}{\wh{\wh \L}}

\newcommand{\arrh}[3]
{
\xymatrix{
{#1} \ar[r]^<<<<{#2}  &{#3}
}
}

\newcommand{\arlh}[3]
{
\xymatrix{
{#1} & \ar[l]_<<<<{#2}  {#3}
}
}

\newcommand{\arrr}[1]
{\arrh {}{#1}{}}

\newcommand{\arr}
{\arrr {}}

\newcommand{\arrl}[1]
{\arlh {}{#1}{}}

\newcommand{\arl}
{\arrl {}}

\newcommand{\arrto}
{\xymatrix{{} \ar@{|-{>}}[r]  & {} } }

\newcommand{\arrinto}
{\xymatrix{{} \ar@{^{(}->}[r]  & {} } }

\renewcommand{\lg}{\wh\L_\Gamma}

\newcommand{\spn}{ S.P. Novikov}

\newcommand{\tra}{\twoheadrightarrow}

\newcommand{\mo}{monomorphism}
\newcommand{\epi}{epimorphism}

\begin{document}

\title
[On the conical Novikov homology]
{On the conical Novikov homology}
\author{Andrei Pajitnov}
\address{Laboratoire Math\'ematiques Jean Leray 
UMR 6629,
Universit\'e de Nantes,
Facult\'e des Sciences,
2, rue de la Houssini\`ere,
44072, Nantes, Cedex}                    
\email{andrei.pajitnov@univ-nantes.fr}

\thanks{} 
\begin{abstract}
Let $\o$ be a Morse form on a closed connected \ma~ $M$.
Let $p:\wh M\to M$ be a regular covering with structure group $G$,
such that $p^*([\o])=0$. The period \ho~
$\pi_1(M)\to\rr$ corresponding to $\o$ factors through a 
\ho~ $\xi:G\to\rr$.
The rank of $\Im \xi$ is called the irrationality degree of $\xi$.
Denote by $\L$ the group ring $\zz G$ and let $\Lxi$ be its Novikov completion.
Choose a transverse $\o$-gradient $v$. 
The classical construction of counting the flow lines of $v$
defines the Novikov complex 
$\NN_*$ freely generated over $\Lxi$ by the set 
%$Z(\o)$ 
of zeroes of $\o$.

In this paper we introduce a refinement of this construction.
We define a subring $\lg$ 
of $\Lxi$
(depending on an auxiliary parameter $\G$ which is a certain cone
in the vector space $H^1(G,\rr)$)
and show that the Novikov complex $\NN_*$ 
is defined actually over $\lg$ and computes 
the homology of the chain complex
$ C_*(\wh M)\tens{\L}\lg   $.
In the particular case when 
$G\approx \zz^2$, and the irrationality degree of $\xi$
equals 2, the ring $\lg$ is isomorphic to the ring 
of  series in $2$ variables $x, y$ 
of the form $\sum _{r\in\nn} a_r x^{n_r}y^{m_r}$
where $a_r, n_r, m_r\in\zz$ and 
both $n_r, \ m_r$ converge to $\infty$ when $r\to \infty$.

The algebraic part of the proof is based on 
a suitable generalization of the classical  algorithm of  approximating 
irrational numbers by rationals.
The geometric part is a straightforward generalization
of the author's proof of the particular case of this theorem 
concerning the circle-valued Morse maps \cite{P-o}.
As a byproduct we obtain a simple proof of the properties of the Novikov
complex for the case of Morse forms of irrationality degree $>1$.

The paper contains two appendices.
In  Appendix 1 we give an overview of the E. Pitcher's 
work on circle-valued Morse theory (1939).
We show that Pitcher's lower bounds for the number of critical points
of a circle-valued Morse map coincide with the torsion-free part of the Novikov 
inequalities (1982). In  Appendix 2 we construct a circle-valued 
Morse map and its gradient \sut~  its unique Novikov incidence 
coefficient is a power series in one variable with an arbitrarily
small convergence radius.

\end{abstract}
\dedicatory{Dedicated to the memory of Andrew RANICKI.}

\maketitle
\tableofcontents

\newcommand{\Lga}{\wh\L_\G}
\newcommand{\Lgam}{{\wh\L}_\G^-}
\newcommand{\Lganul}{\wh\L^\circ_\G}

\newcommand{\smm}{\SS_*(\wh M)}
\newcommand{\sxi}{\wh\SS_*(\wh M, \xi)}
\newcommand{\sxii}{\wh\SS_*(\wh M, \xi_i)}
\newcommand{\sxim}{\wh\SS^-_*(\wh M, F)}
\newcommand{\sximi}{\wh\SS^-_*(\wh M, F_i)}
\newcommand{\sxinul}{\wh\SS^\circ_*(\wh M, \xi)}

\newcommand{\sxig}{\wh\SS_*(\wh M, \G)}
\newcommand{\sximg}{\wh\SS^-_*(\wh M, \FF)}
\newcommand{\sxignul}{\wh\SS^\circ_*(\wh M, \G)}

\newcommand{\cmmm}{\CC_*(\wh M)}
\newcommand{\cxi}{\wh\CC_*(\wh M, \xi)}
\newcommand{\cxii}{\wh\CC_*(\wh M, \xi_i)}
\newcommand{\cxim}{\wh\CC^-_*(\wh M, F)}
\newcommand{\cximi}{\wh\CC^-_*(\wh M, F_i)}
\newcommand{\cxinul}{\wh\CC^\circ_*(\wh M, \xi)}
\newcommand{\cxig}{\wh\CC_*(\wh M, \G)}
\newcommand{\cximg}{\wh\CC^-_*(\wh M, \FF)}
\newcommand{\cxignul}{\wh\CC^\circ_*(\wh M, \G)}

\newcommand{\hol}{ homology }

\newcommand{\steg}{ \SS_*(\wh M)\tens{\L}\Lga }
\newcommand{\cteg}{ \CC_*(\wh M)\tens{\L}\Lga }

\newcommand{\Lximm}{{\wh \L}^-_\xi}
\newcommand{\Lxiimin}{{\wh \L}^-_{\xi_i}}
\newcommand{\Lxii}{{\wh \L}_{\xi_i}}
\newcommand{\Lii}{{\wh \L}_{i}}
\newcommand{\Liimin}{{\wh \L}^-_{i}}

\newcommand{\lili}{\wh\L}
\newcommand{\lilim}{\wh\L^-}
\newcommand{\lilin}{\L^{(n)}}
\newcommand{\nin}{\NN^-_*}

\newcommand{\lin}{\L^-_{i,n}}
\newcommand{\linn}{\L^{(n)}_{i}}
\newcommand{\linol}{\L^-_i}
\newcommand{\lign}{\L^-_{\G,n}}
\newcommand{\lgamin}{\L^-_{\G}}
\newcommand{\lgann}{\L^{(n)}_{\G}}

\section{Introduction}
\label{s:intro}

\subsection{Background}
\lb{su:back}

Apparently the first recorded work  about circle-valued Morse functions 
was Everett Pitcher's paper \cite{Pitcher}, published in 1939.
\footnote{ I am grateful to Andrew Ranicki for pointing out this paper to me.} 
E. Pitcher gave there lower bounds for the number of critical points of a 
Morse map $f:M\to S^1$ 
in terms of homology invariants of the corresponding infinite cyclic covering.
His paper is much less known than  more recent works of S. P. Novikov,
so we outline E. Pitcher's work in Appendix 1.
We show in particular that his inequalities are equivalent to 
the torsion-free part of the Novikov inequalities.

The circle-valued Morse theory as it exists today was initiated by \spn~
 \cite{NovikovDAN}. In this work  \spn~ considered a Morse map $f:M\to S^1$
and outlined a construction of a chain complex generated 
by critical points of $f$. The base ring of this complex is the ring 
of integer Laurent series with finite negative part
$\zz((t))=\zz[[t]][t^{-1}]$. The chain complex computes the {\it semi-open}
homology $H_*(\wh M, \infty_+)$
of the corresponding infinite cyclic covering $\wh M$.
%\footnote{See \S ??? for  a discussion of this and similar notions.}
One immediate consequence is the Novikov inequality
for the number of critical points of $f$ of index $r$:
\bq\lb{f:nov_ineq}
m_r(f)\geq 
b_r(M, \xi)+
q_r(M, \xi)+
q_{r-1}(M, \xi),
\end{equation}
where we denote by $b_r$ and $q_r$  the rank and torsion number of
the above homology in degree $r$, and   $\xi$ is  the homotopy class of $f$
in the group $[M, S^1]\approx H^1(M,\zz)$.
\spn~ suggested also a generalization
of this construction to the case of Morse forms
(recall that a closed 1-form is called {\it Morse form}, if locally it is the differential 
of a Morse function).
The ring $\zz((t))$ is replaced in this case
by a corresponding completion of the ring $\zz[\zz^k]$,
where $k$ is the maximal number of 
rationally independent periods of the Morse form in question
(i.e. the irrationality degree 
of the form).

The details of the construction of the chain complex above
were not clarified in \cite{NovikovDAN}. Several authors suggested different approaches 
to construction of this chain complex and proof of its properties.
In an influential work \cite{FarberExN}
M. Farber proved the exactness of the inequalities
\rrf{f:nov_ineq}
for the case of manifolds of dimension $\geq 6$ with fundamental group $\approx \zz$.
He also gave a proof of the Novikov inequalities without using the properties
of the chain complex above. Instead of using the Laurent series ring,
M. Farber works with a suitable localization
of the polynomial ring $\zz[t]$; the numerical invariants
derived from the localized homology of the infinite cyclic covering
coincide with the Novikov numbers.

In his PhD thesis \cite{SikoravThese} J.-Cl. Sikorav 
suggested a generalization of the Novikov homology,
replacing the free abelian covering 
of the manifold in question by the universal covering.
This homology (called sometimes {\it Novikov-Sikorav homology} )
turned out to be sufficiently strong 
so as to detect    3-manifolds fibered over a circle.

A construction of the Novikov complex for a circle-valued Morse map
based on the counting of the flow lines of the gradient of the Morse map 
was given in the author's work \cite{P-o}.
This method was inspired by E. Witten's framework for Morse theory
\cite{Witten1982}. In the work \cite{P-o} the first complete proof of the properties
of the Novikov complex was given. The proof uses the classical 
Morse theory on compact submanifolds of the infinite cyclic covering,
and inverse limit arguments. 

A different approach was developed in the articles of M. Pozniak
%\cite{Pozniak1}, 
%\cite{Pozniak2} 
(PhD thesis at the University of Warwick, 1994, published later in 
\cite{Pozniak2}) 
and F. Latour \cite{Latour}.
Both use the Floer's method \cite{Floer}.
The proof splits into two parts. Firstly one proves that
$d^2=0$ by studying the space of broken flow lines between critical points. 
M. Pozniak's argument \cite{Pozniak2} generalize the argument of D. Salamon \cite{Salamon}
concerning the real-valued Morse functions. It is based on the theory of Fredholm
operators in Sobolev spaces. F. Latour's argument \cite{Latour}, \S 2
is based on the transversality properties of 
manifolds with singular  boundary.
The second step is the computation of the homology of the resulting complex.
It is done by different versions of Floer's continuation method.
One shows that the homology of the Novikov complex does not change 
while deforming the Morse form within its cohomology class.
Afterwards one shows that in a given de Rham cohomology class 
there exists always a Morse form having the same 
Novikov complex as the completed Morse complex of  a real-valued Morse function.
This second step was carried out by M. Pozniak in his thesis \cite{Pozniak2},
and by F. Latour in the Chapter 2 of \cite{Latour}.

\subsection{Overview of the paper}
\lb{su:over}

The main object of study in this paper is the Novikov complex
for  Morse forms of irrationality degree $k>1$. 
Let $\o$ be such form, and $p:\wh M\to M$  a regular covering 
with structure group $G$, \sut~ $p^*(\o)$ is cohomologous to zero.
Let $v$ be a transverse 
$\o$-gradient. Our approach to construction of the Novikov complex
is based on approximation of $\o$ by rational Morse forms
(that is, multiples of differentials of circle-valued Morse maps).
Choose  rational Morse forms $\o_1, \ldots, \o_k$ 
approximating $\o$, so that $v$ be an $\o_i$-gradient for each $i$.
The Novikov complex for $\o$ is then defined 
over a ring that we call {\it the conical completion } of the group ring $\zz G$.
\footnote{These completions were present implicitly already  in the author's paper 
\cite{Pepr2}.}
This conical completion is a much smaller ring than the initial Novikov ring.
For a case when the homology classes 
%$\xi_1, \ldots, \xi_k$ 
of the forms  $\o_1, \ldots, \o_k$ form a {\it regular family} (see Sections \ref{s:reg-cc},
\ref{s:twisted} for the definition)
this ring is a subring of the ring of twisted Laurent series in $k$ variables.
This subring consists of Laurent series that we call {\it special 
Laurent series} (see Section \ref{s:twisted} ). 

Theorem \ref{t:exist-regular} (Section \ref{s:reg-cc})
says that we can always find rational forms $\o_1, \ldots, \o_k$
arbitrarily close to $\o$,
\sut~ the family of their homology classes
is regular. The argument is based on a suitable generalization
of the classical algorithm of approximating irrational
numbers by rationals. 

The properties of the conical refinement
of the Novikov complex are stated in the main theorem of the paper
(Theorem \ref{t:main_th}); the proof is given in Section \ref{s:prf_main}.
The advantage of our approach  is that 
the proof of the properties of the Novikov complex is reduced 
to the proof for the case of rational forms 
done in \cite{P-o}. 
The section \ref{s:circle}
contains the outline of this proof.

The paper contains two appendices.
In  Appendix 1 we give an overview of the E. Pitcher's 
work on circle-valued Morse theory (1939).
We show that Pitcher's lower bounds for the number of critical points
of a circle-valued Morse map coincide with the torsion-free part of the Novikov 
inequalities (1982). In  Appendix 2 we construct a circle-valued Morse map and 
its gradient \sut~ its unique Novikov incidence coefficient 
is a power series in one variable with an arbitrarily small convergence radius. 

\subsection{Remarks on terminology}
\lb{su:rem-term}

\noi 
In this paper we work with three types of regular coverings 
of the manifold $M$:
1) infinite cyclic coverings,
2) covering with free abelian structure group,
3) coverings with a structure group $G$ endowed with an epimorphism onto a free abelian group.
The coverings of the first type will be denoted by $\ove M\to M$
(or  $\ove M_i \to M$ if we work with several such coverings).
The coverings of the second type will be denoted by $\wi M\to M$.
The coverings of the third  type will be denoted by $\wh M\to M$.
We assume that the structure groups act on regular coverings from the right,
so that the singular chain complex of the covering is a  right module 
over the group ring.
The singular chain complex of a topological space $X$ is denoted by 
$\SS_*(X)$, the simplicial chain complex of a simplicial space $Y$
is denoted by 
$\D_*(Y)$.
%\pa
%\noi 
The end of a proof is marked by $\square$,
the end of a remark or a definition is marked by $\triangle$.

\subsection{Acknowledgements}
\lb{su:ack}

I am indebted  to A. Ranicki for many discussions
on  circle-valued Morse theory.
I am grateful to G. Ziegler  for the references about integral cones 
and to J. Gubeladze for nice and helpful discussion about bases in integral cones.
Thanks to 
anonymous referee,  whose remarks have lead to a 
considerable improvement of the manuscript.
Many thanks to F. Bogomolov for his constant support.

\section{Conical completions of group rings and chain complexes}
\label{s:cc}

\subsection{Completions of group rings}
\lb{su:c-rings}

Let $G$ be a group; we denote by $\L$ the group ring $\zz G$.
Let $\xi:G\to\rr$ be a \ho.
Intuitively, the Novikov completion $\Lxi$
of the ring $\L$ consists of some special infinite linear combinations
of the elements of $G$, namely the combinations that are infinite
in the direction of decreasing of $\xi$. To give a precise definition
let $\lLL$ be the set of all formal linear combinations
(infinite in general) 
$\l=\sum_{g\in G} n_g g, \ n_g\in \zz$.
For $\l\in \lLL$ put $\supp \l = \{g~|~ n_g\not=0\}$.
For $C\in\rr$ put 
$$
[\l]_{\xi, C}
=
\{g\in \supp \l~|~  \xi(g)\geq C\}.$$

\bede\lb{d:l-xi}(\cite{NovikovDAN}, \cite{SikoravThese})
\bq\lb{f:lxi}
\Lxi=\{\l\in\lLL~|~ \forall C {\rm \ \ the \ \ set \ \ } [\l]_{\xi, C} {\rm \ \ is \ \ finite \ \ }\}
\end{equation}
\bq\lb{f:lhat-min}
\lLL_\xi^-=\{\l\in\lLL~|~ \supp \l \sbs \xi^{-1}(]-\infty, 0]) \}
\end{equation}
\bq\lb{f:lxim}
\Lximm = \Lxi \cap \lLL_\xi^-.
%\Lximm=\{\l\in\lLL~|~ \forall C {\rm \ \ the \ \ set \ \ } [\l]_{\xi, C} {\rm \ \ is \ \ finite \ \ }\}
\end{equation}

\end{defi} 
Equivalently
\bq\lb{f:lxi-bis}
\Lxi=\left\{\l= \sum_{i=0}^\infty n_i g_i~|~n_i\in\zz, \ \ \xi(g_i)\to -\infty  {\rm \ as  \ \ } i\to\infty \right\}
\end{equation}
\bq\lb{f:lxim-bis}
\Lximm=\left\{\l= \sum_{i=0}^\infty n_i g_i~\big|~n_i\in\zz,\ \xi(g_i)\leq 0   
{\rm \ for \ all  \ \ } i, {\rm \ and  \ \ }  \ \xi(g_i)\to -\infty  {\rm \ as  \ \ } i\to\infty  \right\}
%\Lximm=\{\l\in\lLL~|~ \xi|\supp\l\leq 0 {\rm \ \ and  \ \ }\l= \sum_{i=0}^\infty n_i g_i~|~n_i\in\zz, \ \ \xi(g_i)\leq 0, \ \ \xi(g_i)\to -\infty    \}
\end{equation}

It is easy to see that $\Lxi$ is a ring, and $\Lxim$ is a  subring of $\Lxi$.
These rings will be called {\it Novikov completions} 
of the group ring $\L$. $\qt$

\pa\noindent
The basic algebraic objects of the present paper 
are  {\it conical Novikov completions } of the group rings,
introduced in the next definition.
%\footnote{These completions were implicit already  in the author's paper 
%\cite{P-epr2}.}

\bede\lb{f:lgas}
Let $\G=\{\xi_1, \ldots, \xi_k\}$ be a finite set of \ho s $G\to\rr$.
Put 
\begin{gather}\lb{f:def-lgas}
\Lga=\bigcap_i\wh\L_{\xi_i}, \ \ \Lgam = \bigcap_i {\wh\L}^-_{\xi_i}, \\
\Lganul = \bigcup_{g\in G}g\cdot \Lgam. \ \ \ \ \ \qt
%\Lganul = \{\l\in\lLL ~|~ \exists  \ g\in G {\rm \ \ with \ \ } g\l\in \Lgam \}.
\end{gather}

\end{defi}

\bere\lb{r:defs-limits}
The definitions of the rings $\Lximm$ and $\Lgam$
can be reformulated in terms of inverse limits of rings.
Namely, for $\l\in \L$ define the {\it $\xi$-height } $h_\xi(\l)$ of $\l$ as follows:
$$
h_\xi(\l) = \max_{g\in \supp\l} (\xi(g)).
$$
For a real number $K$ put $\L^-_{\xi, K} = \{\l\in \L~|~ h_\xi(\l) \leq K\}$.
We abbreviate $\L^-_{\xi, 0}$ to $\L^-_{\xi}$. 
Then $\Lximm$ is isomorphic to the inverse limit of the following sequence of ring projections
$$
%\L^-_{\xi} \arl 
\L^-_{\xi}\Big/\L^-_{\xi, -1}\arl \ldots \arl\L^-_{\xi}\Big/\L^-_{\xi, -n} \arl\ldots 
$$
For  a finite set $\G=\{\xi_1, \ldots, \xi_k\}$ of \ho s $G\to\rr$
put $\L^-_{\G, K} = \bigcap_i \L^-_{\xi_i, K}$.
We abbreviate $\L_{\G, 0}$ to $\L^-_{\G}$. 
Then $\Lgam$ is isomorphic to the inverse limit of the following sequence of ring projections
\begin{flalign*}
 &&
 %\L^-_{\G} \arl 
 \L^-_{\G}\Big/\L^-_{\G, -1}\arl \ldots \arl
 \L^-_{\G}\Big/\L^-_{\G, -n} \arl\ldots 
  && \qt
 \end{flalign*}

\enre

We have $\Lgam\sbs \Lganul\sbs \Lga$.
Assume that all $\xi_i$ factor through 
an epimorphism 
$\z:G\to \zz^n$,
that is, there are $\bar\xi_i:\zz^n\to\rr$, \sut~
$\xi_i=\bar\xi_i\circ \z$
(this is always the case for some $n$, if $G$ is finitely generated).
Choose and fix such an epimorphism $\z$. 
%For a given finite set $\G=\{\xi_1, \ldots, \xi_k\}$ as above
%there exists always an epimorphism $p:G\to \zz^k$
%\sut~ all $\xi_i$
%factor through $p$, that is, there are $\bar\xi_i:\zz^k\to\rr$, \sut~
%$\xi=\bar\xi_i\circ p$.
%Choose and fix such an epimorphism $\z$. 
The extensions of $\bar\xi_i$ to linear forms $\rr^n\to\rr$
will be denoted by the same symbols $\bar\xi_i$.
Put
\bq\lb{f:c-gamma}
C_\G=\left\{x\in \rr^n~|~ \bar\xi_i(x)\leq 0
\ \ {\rm for\ \  every } \ \ i
\right\}.
\end{equation}

\bepr\lb{p:cones}
Assume that $C_\G$ is a solid cone. 
\footnote{ Recall \cite{BoydVandenConv} 
that a subset $X\sbs \rr^n$ is called {\it cone} if 
for every $a\in X$ and $\t\geq 0$ we have $\t a\in X$. 
A cone is called {\it solid} if it has a non-empty interior. }
Then $\Lganul=\Lga$.
\enpr
\Prf
Let $\l\in\Lga$. The set $E=\supp\l \ \sm  \ \zeta^{-1}(C_\G)$
is finite. Since $C_\G$ is a solid cone, there is a ball 
$B\sbs \Int(C_\G)$ of any given radius. 
So we can assume that for some $a\in \z(G)$ 
we have $\z(E)+a \sbs \Int(C_\G)$. Therefore the support of
$\z^{-1}(a)\cdot E$ is in $C_\G$, and the proposition is proved. $\qs$

\bede\lb{d:admiss}
We say that $\G$ is {\it admissible}, if  $C_\G$ is a solid cone.
\end{defi}

\bere\lb{re:re}
Observe that $C_\G$ is a solid cone if and only if
there exists $g\in G$ \sut~ $\xi_i(g)<0$ for every $i$.
Thus for a family $\G$ the property of being admissible
does not depend on the choice  of the epimorphism $\zeta$.
$\qt$
\enre

A basic example of an admissible cone is provided 
by the following obvious proposition.
\bepr\lb{p:admiss}
Let $\xi:G\to\rr$ be a \ho, that factors through
$\z:G\to\zz^n$. Then there is
a \nei~ $\UU$ of $\bar\xi$ in $L(\rr^n, \rr)$
\sut~ any set 
$\G=\{\xi_1, \ldots, \xi_k\}$
with all $\bar\xi_i$ in $\UU$ is admissible.
If $||\cdot ||$ 
denotes the norm in 
$L(\rr^n, \rr)$
associated to a scalar product in $\rr^n$,
then we can take \
$\UU=\{\eta~|~ ||\eta-\bar\xi||<||\bar\xi||\}$. $\qs$
\enpr

\bere\lb{r:nonadm}
For  non-admissible cones the conclusion of 
Proposition \ref{p:cones} is not valid in general.
For example, let $G=\zz,\ \  \xi=\id:G\to\zz$,
and $\G=\{\xi, -\xi\}$. Then $\Lga=\zz[\zz],\ \  \Lgam\approx\zz\ \ $,
and $\Lganul $ is not a subgroup. $\qt$
\enre 

\subsection{Completions of chain complexes}
\lb{su:c-compl}

Similarly to  completions of group rings we 
can define completions of singular and cellular 
chain complexes of coverings.
Let $M$ be a connected topological space, and 
$p:\wh M\to M$ a regular covering of $M$
with structure group $G$.
Let $\SS_*(\wh M)$ be the singular chain complex of $\wh M$.
We denote by $\wh{\wh{\SS}}_*(\wh M)$ 
the set of all formal linear combinations
(infinite in general) 
$T=\sum_{\s} n_\s \s, \ n_\s\in \zz $
where $\s$ ranges over singular simplices of $\wh M$. 

Let $\xi:G\to\rr$ be a \ho~
factoring through some epimorphism
$G\to\zz^n$. We have then a 
regular covering $\wi p: \wi M \to M$ with 
structure group $\zz^n$, and a 
commutative diagram of $\zz^n$-coverings 
$$
\xymatrix{ 
\wi M  \ar[d]^{\wi p}  \ar[r]^{\wi \r} & \rr^n \ar[d]^{p_0}\\
M \ar[r]^\r & \ttt^n
}
$$
The \ho~ $\xi$ can be considered as an element of $H^1(M,\rr)$. 
We have $\xi=\r^*(\xi_0)$ for some class 
$\xi_0\in H^1(\ttt^n,\rr)$.
Let $\o_0\in\O^1(\ttt^n)$ be a closed 1-form
with $[\o_0]=\xi_0$, choose a function $\phi:\rr^n\to\rr$ 
\sut~ $d\phi=p_0^*(\o_0)$, and let $F=\phi\circ\wi\r$.
Observe that $F(gx)=F(x)+\xi(g)$, where $g\in G$.
%For $C\in\rr$ 
%put
%$[T]_{F,C}=\{\s\in\supp T~|~ (F|\s) \geq C\}$.

\bede\lb{d:sing-c}
(\cite{NovikovDAN}, \cite{SikoravThese})
For $C\in\rr$ and 
$T=\sum_{\s} n_\s \s \in \wh{\wh{\SS}}_*(\wh M)$
put $\supp T = \{\s\in \SS_*(\wh M)~|~ n_\s\not=0\}$.
Let 
$[T]_{F,C}=\{\s\in\supp T~|~ (F|\s) \geq C\}$. Put

\bq\lb{f:sxi}
\sxi=\left\{T\in\wh{\wh{\SS}}_*(\wh M)~|~ 
\forall C {\rm \ \ the \ \ set \ \ } [T]_{F, C} {\rm \ \ is \ \ finite \ \ }\right\}
\end{equation}
\bq\lb{f:shat-min}
\wh{\wh{\SS}}^-_*(\wh M, F ) =
\left\{T\in \wh{\wh{\SS}}_*(\wh M) ~|~
\supp T \sbs F^{-1}(]-\infty, 0]) \right\}
\end{equation}
\bq\lb{f:sxim}
\sxim=\sxi \ \cap \ \wh{\wh{\SS}}^-_*(\wh M, F)
\end{equation}
These chain complexes will be called {\it Novikov completions}  
of the singular chain complex of $\wh M$.
$\qt$
\end{defi}

Observe that the completion \rrf{f:sxi} depends only on 
$\xi$ but not on the particular choices of 
$\o_0$ in the de Rham cohomology class $\xi_0$, 
neither on the choice of the function 
$\phi$.

\bede\lb{f:chain-cones-cc}
Let $\xi_1, \ldots, \xi_k : G\to\rr$
be \ho s all factoring through
$\z:G\to\zz^n$. For each $i$ choose the corresponding function
$\phi_i:\rr^n\to\rr$, and put $F_i=\phi_i\circ \wi\rho$.
%Let $\o_1, \ldots, \o_n\in\O^1(M)$ be  closed 1-forms,
%\sut~ 
%$p^*(\o_i)=dF_i$,
%where $F_i:\wh M\to\rr$ are $\smo$ functions.
%Put $\xi_i=[\o_i]$,
%then $\xi_i$ can be considered as a \ho~ $G\to\rr$.
Write $\G=\{\xi_1,\ldots, \xi_k\}$, 
and $\FF=\{F_1,\ldots, F_k\}$.
Put

\begin{gather}\lb{f:sxig}
\sxig=
\bigcap_i \ 
\sxii, \ \ \ \ \ \ 
%\lb{f:sximg}
\sximg=
\bigcap_i \ 
\sximi, \\
%\lb{f:sxignul}
\sxignul=
\bigcup_{g\in G} g\cdot \sximg.
\end{gather}

\end{defi}
We have obviously $\sximg \sbs \sxignul \sbs \sxi$. $\qt$

\bere\lb{r:inv-lim-def}
For $r\in\nn$ 
let $\wh M^{(r)} = \{x\in \wh M ~|~ \forall i \ \ F_i(x)\leq -r \}$.
We obtain a decreasing filtration 
$\FF_r= \SS_*(\wh M^{(r)})$ in $\SS_*(\wh M)$
and $\sxig$ is the inverse limit of the corresponding 
inverse system
$$
\ldots \arl \SS_* \big(\wh M\big)/ \SS_* \big(\wh M^{(r)}\big) 
\arl 
\SS_* \big(\wh M\big)/ \SS_* \big(\wh M^{(r+1)}\big)  \arl \ldots 
\hspace{1cm}\qt 
$$  
\enre

Observe that $\sximg$ is a right $\Lgam$-module.
The natural inclusion 
$\sximg$ $ \arrinto  \sxig$ 
extends obviously to a \ho~ 
\bq\lb{f:iii}
I:
\sximg \tens{\Lgam}\Lga
\arr
\sxig.
\end{equation}

%The proof of the next proposition is similar to the proof of 
%Proposition \ref{p:cones} and will be omitted.
\bepr\lb{p:cones-tensor}

Assume that $C_\G$ is a solid cone. Then 

\been\item
$\sxignul =  \sxig$.
\item 
%The \ho~ 
$I$ is an isomorphism.
\enen
\enpr 
\Prf
The proof of the first part and of surjectivity of $I$ 
is similar to Proposition \ref{p:cones} and will 
be omitted. As for the injectivity of $I$, 
assume that for some non-zero element 
$\xi=\sum_{j=1}^N g_j\otimes\xi_j$ 
(with $g_j\in G, \ \xi_j\in \sximg$ )
we have $I(\xi)=0$.
Pick $g\in G$ \sut~ $gg_j\in \lgamin$ for every $j$.
Then $g\xi = \sum_j gg_j\otimes \xi_j \in \sximg$,
and $I(g\xi)=gI(\xi)=0$,
which is impossible since $I~|~\sximg$ is 
injective. $\qs$

\subsection{Singular homology versus cellular homology}
\lb{su:sing-vers-cell}

In this subsection we will use the terminology of the previous one. 
If $M$ is a CW complex, then we can endow 
$\wh M$ with a $G$-invariant CW structure. 
Similarly to Definition \ref{d:sing-c} and Remark \ref{r:inv-lim-def}
one can define conical completions of the 
cellular chain complex $\CC_*(\wh M)$:
\bede\lb{d:cell-compl}
Let  $\CC^{(r)}_*(\wh M)$ be the 
subcomplex of all cellular chains
contained in 
$\wh M^{(r)}$.
These subcomplexes form a decreasing filtration
of $\CC_*(\wh M)$; denote by $\wh\CC_*(\wh M, \G)$
the corresponding inverse limit.
The subcomplexes $\CC^{(r)}_*(\wh M)$ for $r\geq 0$
form a decreasing filtration
of $\CC^{(0)}_*(\wh M)$; denote by $\cximg$ 
the corresponding inverse limit.
Put 
$$
\cxignul
=
\bigcup_{g\in G} g\cdot \cximg
\sbs 
\cxig. \hspace{1cm}\qt 
$$ 
\end{defi}
We are going to prove that the homology modules of the cellular and singular 
versions of the conical completions  are isomorphic. 
Let $M^{[m]}$ be the $m$-skeleton of $M$.
We have canonical $\zz G$-equivalences:
$
\Phi: \CC_*(\wh M) \arr \SS_*(\wh M) 
$
\sut~ 
\been\item 
$\Phi(\CC_m(\wh M))\sbs \SS_*(\wh M^{[m]})$;
\item
if we endow $\CC_*(\wh M)$ with the filtration by
the subcomplexes 
$$\CC_*^{(s)}\ = \ \{0\arl \CC_0(\wh M)\arl \ldots \arl \CC_s(\wh M)\arl 0 \ldots   \}$$
(we call this filtration {\it trivial} ),
then the map $\Phi $ induces the identity 
\iso~ in the homology of the quotient complexes:
$$
\CC_m(\wh M)
\arrr \approx 
H_m(\wh M^{[m]}/\wh M^{[m-1]}).
$$
\enen 

(see for example \cite{P-o}, Section 3).
Extending  $\Phi$  to the completion
we obtain a chain map
\begin{equation}
\lb{f:mapss}
\wh\Phi:
\cxig
\arr
\sxig.
\end{equation}
%\ \ \ \ 
%\wh\Phi^-:
%\cximg
%\arr
%\sximg, \\ 
%\lb{f:mapss1}
%\wh\Phi^\circ:
%\cxignul
%\arr
%\sxignul
%\end{gather}

\bepr\lb{p:isos}
Assume that $M$ is a finite CW-complex 
and $C_\G$ is a solid cone.
then the chain map \rrf{f:mapss}
induces an \iso~  in homology.
\enpr
\Prf
Both $\cxig$
and
$\sxig$
are inverse limits 
of the quotients 
of their
filtrations.
It will be convenient to replace the filtration in 
the singular chain complex by another one which is 
equivalent to it.
Denote by 
$\wh M_r$ 
the minimal CW subcomplex of $\wh M$ 
containing all the cells in $\wh M^{(r)}$.
We obtain filtrations in $\CC_*(\wh M)$ and 
$\SS_*(\wh M)$
equivalent to the previously introduced ones.
The chain equivalence $\Phi$ gives rise to a map of corresponding
inverse systems which induces an isomorphism of the inverse limits of 
homology modules. 
Both filtrations satisfy obviously the Mittag-Leffler 
condition,
therefore $\lim^1$ vanish for both, and the homology
modules of the inverse limits are isomorphic. $\qs$

\subsection{ Completions versus tensor products}
\lb{su:compl-tens}

Along with conical completions 
\newline
$\sxig$ we can consider tensor products
$ \SS_*(\wh M)\tens{\L}\Lga$.
We have a natural inclusion 
$\steg\arrinto \sxig$.
\bepr\lb{p:compl-tens}
Assume that $M$ has a homotopy type of a finite CW complex.
Then the inclusion 
$i:
\steg
\arrinto
\sxig
$
induces an \iso~ in \hol.
\end{prop}
We can assume that $M$ is a finite CW complex.
Observe that (contrarily to the case of singular homology)
the inclusion $\cteg\arrinto\cxig$ is an \iso.
Consider the following diagram
$$
\xymatrix{ 
 {\steg}  \ar[r]^{i} & \sxig \\
\cteg  \ar[u]^{\Phi\otimes \Id}  \ar[r]_\approx & \cxig   \ar[u]^{\wh\Phi}
}
$$
Both vertical arrows are induced by the canonical  $\zz G$-equivalence 
$\Phi: \CC_*(\wh M)\sim S_*(\wh M)$.
 The map $\wh\Phi$ is  a homology equivalence by 
 Proposition
 \ref{p:isos},
 therefore $i$ is also a homology equivalence.
 $\qs$

For infinite CW complexes the inclusion $i$ above 
is not necessarily a chain equivalence.
\begin{exam}\lb{e:infin}
Let $X$ be an infinite wedge of 2-spheres and $M=S^1\vee X$.
Let $\xi$ be the generator of $H^1(M)$
and $\G=\{\xi\}$. Then the map
$ H_2(\wh M)\tens{\L}\Lga \arrto 
H_2(\sxig)$
is not epimorphic. $\qt$ 
\end{exam}

\bere\lb{r:histo-tens}
In the original article \cite{NovikovDAN}
\spn~ works with completions of singular chain complexes. 
The same approach is used by J.-Cl. Sikorav in 
\cite{SikoravThese}.
In several later papers 
the authors use another version of the Novikov homology, 
defined via tensor products. The above proposition imply that 
both constructions are equivalent in the case of finite CW complexes. 
$\qt$
\enre

%\newpage
\section{Regular conical completions}
\label{s:reg-cc}

In our applications we will be working with some special families 
of \ho s 
$G\to \rr$. For these families the corresponding  conical completions
are isomorphic to certain subrings of  
 the ring of power series in several variables 
(non-commuting in general).

\subsection{$\xi$-regular families}
\lb{su:xi-reg}

\bede\lb{d:xi-reg}
Let 
$\xi:\zz^k\to\rr$ 
and $\xi_1, \ldots, \xi_k : \zz^k\to\zz$ be group \ho s.
We say that the family 
$\G=\{\xi_1, \ldots, \xi_k\}$
is $\xi$-regular, if
\been\item
$
(\xi_1, \ldots, \xi_k)
$
is a basis in the free abelian group $\Hom(\zz^k,\zz)$.
\item 
All the coordinates of $\xi$ in the basis 
$
(\xi_1, \ldots, \xi_k)
$
are strictly positive.  
\enen 
\end{defi}
It is easy to see that for a $\xi$-regular
family
$\G$ the cone $C_\G$ is a  solid cone,
and 
$C_\G\sbs \xi^{-1}\big(]-\infty, 0]\big)$.
For a vector $v$ in a vector space we denote by $l(v)$
the ray generated by $v$,
that is, $l(v)=\{tv~|~ t\geq 0\}$. 
For two rays $l_1, l_2$ in the Euclidean vector space 
$\rr^k$ we denote by $d(l_1, l_2)$ the usual 
angular distance between $l_1$ and $l_2$
(recall that $d(l_1, l_2)$ is defined as the 
distance between the intersection points $l_1\cap S$
and $l_2\cap S$ where $S$ is the Euclidean sphere centered in $0$ of radius $1$,
and the distance is with respect to the induced Riemannian metric
on $S$).

\beth\lb{t:exist-regular}
Let $\xi:\zz^k\to\rr$ be a monomorphism, and $\e>0$.
Then there is a $\xi$-regular family 
$\G=\{\xi_1, \ldots, \xi_k\}$
\sut~ $d(l(\xi),l(\xi_i))<\e$ for every $i$. 
\enth 

Choosing a $\zz$-basis in the group 
$\zz^k$ we identify the group $\Hom(\zz^k,\zz) $
with $\zz^k$; the vector space $L(\rr^k, \rr)$ 
is then identified with $\rr^k$. 
The theorem above  follows from the next theorem 
dealing with vectors in $\rr^k$.

\beth\lb{t:exist-regular-bis}
Let $v$ be a vector in $\rr^k=\zz^k\otimes \rr$
\sut~  its coordinates 
$v^{(1)}, \ldots , v^{(k)}$ 
are linearly independent over $\qq$.
\footnote{ Such vectors will be called {\it maximally irrational}.}
Assume that $v^{(1)}>0$. Let $\e>0$. 
Then there are vectors 
$u_1, \ldots, u_k\in \zz^k$ 
\sut~ 
\been\item The family $\FF=(u_1, \ldots, u_k) \in\zz^k$
is a $\zz$-basis of \  $\zz^k$.
\item
$
v=\alpha_1u_1 + \ldots +\alpha_k u_k,
$
with $\alpha_i>0$ for every $i$,
\item 
$d(l(v), l(u_i)) <\e$ for every $i$,
\item the first coordinates of all vectors $u_i$ are $>0$.
\enen
\enth 
\Prf
Let us start by constructing a family $\GG$
satisfying the properties 2) -- 4) of Theorem 
\ref{t:exist-regular-bis}.
Let $N\in \nn$ be a natural number with 
$N>\frac {\sqrt{k-1}}{\e}$. 
Denote by $H$ the hyperplane 
$\{x~|~x^{(1)}=N\}$.
Pick $k$ vectors 
$a_1, \ldots, a_k \in H\cap \zz^k$
\sut~ the points $a_i$ form a $(k-1)$-simplex in $H$
of diameter $\leq \sqrt{k-1}$, containing 
the point $l(v)\cap H$ in its interior. 
It is obvious that 
the family $\GG=(a_1, \ldots, a_k)$
satisfies the properties 2) and 4) of Theorem
\ref{t:exist-regular-bis}.
As for the property 3), recall that
the central projection of a sphere onto its tangent hyperplane
is a length-increasing map, 
therefore for any rays $l_1, l_2$ having non-empty 
intersection with $H$
we have
$d(l_1, l_2) 
\leq \frac {\sqrt{k-1}}  N ||l_1\cap H -l_2\cap H||$; 
the property 2) follows. 
Proceeding to the property 1),
observe first that $v$ being maximally irrational,
the property 2) guarantees that the family
$\GG$ is a basis in $\rr^k$.
To achieve the property 1) 
we need to refine the above construction.

Recall that a family $\BB$ of vectors in $\zz^k$ is called 
{\it unimodular} if it is a $\zz^k$-basis.
The cone $C(\BB)$ is called {\it unimodular} if $\BB$ is unimodular.
The theorem \ref{t:exist-regular-bis} 
follows immediately from the following well-known fact
(see the book of W. Bruns and J. Gubeladze 
\cite{BrunsGubeladze}, Th. 2.72 for a proof of a more general result).

\beth\label{t:cone-union-unimod}
For every family $\BB$ of vectors in $\zz^k$ the cone $C(\BB)$ 
is a union of unimodular cones.
$\qs$
\enth 

To make our exposition self-contained we will give  the 
full proof of our theorem
\ref{t:exist-regular-bis}; 
the argument below is essentially equivalent to the proof of 
theorem \ref{t:cone-union-unimod} given in \cite{BrunsGubeladze}.
\footnote{After this article was completed, I became aware that
a similar argument was also used by D. Sch\"utz \cite{S}
in his work about $K$-theory of Novikov rings.}

For a family $\FF=(v_1, \ldots , v_k)$ of vectors in $\rr^k$ 
we denote by $P(\FF)$ the parallelotope generated by $\FF$, that is
$$
P(\FF)=\left\{\sum_i \l_i v_i ~|~  0\leq \l_i \leq 1\right\}.
$$
We denote by $P_0(\FF)$ the {\it semi-open parallelotope} generated by $\FF$, that is
$$
P_0(\FF)=\left\{\sum_i \l_i v_i ~|~ 0\leq \l_i < 1\right\}.
$$
Assume that $\FF\sbs \zz^k$. A parallelotope $P(\FF)$ is called {\it empty}
if it contains no point of the lattice $\zz^k$ except its vertices. 
We will use the following simple lemma (the proof is omitted).

\bele\lb{l:empty-parallelotope}
Let $\FF$ be a family of $k$ linearly independent vectors in $\zz^k$.
The following  properties are equivalent:
\been\item
$P(\FF)$ is empty.
\item $P_0(\FF)\cap \zz^k=\{0\}.$
\item $|\det \FF| = 1.$
\item $vol(P(\FF))=1.$ $\qs$
\enen 
\enle

A family of vectors in $\zz^k$ 
satisfying the conditions 2) -- 4) of
Theorem \ref{t:exist-regular-bis}
will be called {\it admissible}.
In the set of all admissible families choose a family 
$\FF_0=(a_1, \ldots, a_k)$
\sut that 
$|\det(\FF_0)|$
is minimal possible. 
I claim that this volume equals $1$.
Indeed, assume that 
$vol(\D) > 1$.
Then $P_0(\FF)$ contains at least
one point $q\in\zz^k$ different from $0$.
Two possibilities can occur:

1) The point $q$ is in one of the semi-open edges 
$[0, a_1[, \ldots, [0, a_k[$ of $P(\FF_0)$,
say in $[0, a_1[$. In this case
replacing the family 
$\FF_0$ by the family 
$\FF=(q, a_2, \ldots, \a_k)$
diminishes the volume of the corresponding parallelotope. 

2) The vector $q$ is not collinear to any of the vectors 
$a_1, \ldots, a_k$.
Write $q$ as a linear combination
of vectors $a_1, \ldots, a_k$
and let $r$ be the the number of non-zero coefficients in this 
linear combination. We can assume that 
$q=\b_1a_1+\ldots +\b_ra_r$
with $1>\b_i>0$ for every $i$. Observe that $r\geq 2$.
Replacing in the family $\FF_0$ the vector $a_j$ by $q$
(here $j\leq r$)
we obtain a new family $\FF_j$
satisfying obviously the conditions 3) and 4)
of Theorem \ref{t:exist-regular-bis}.

\bele\lb{l:lin-alg}
1)  For every $1\leq j\leq r$ 
the family $\FF_j$ is a basis in $\rr^k$.

2) The union of the cones $C_{\FF_j}$ equals $C_{\FF_0}$.

3) For every $j$ we have 
$|\det(\FF_j)|< 
|\det(\FF_0)|$.
\enle
\Prf
The proof of the 
points 1) and 2) of the
lemma is an easy 
argument from linear algebra and we will omit it.
Let us just outline the geometric contents of the point 2.
Consider the $(k-1)$-simplex $K$ in $H$ with vertices 
$
l(a_1)
\cap H
, \ldots, l(a_k)
\cap H
$. 
Adding one more vertex 
$l(q)\cap H$ we obtain a simplicial subdivision 
of $K$ containing $r$ simplices of dimension $k-1$.

As for the point 3), let us show for example 
that 
$
|\det(\FF_1)|< 
|\det(\FF_0)|$.
We have 
$|\det(q, a_2, \ldots , a_k)|
=
\b_1 
|\det(a_1, a_2,  \ldots , a_k)|
<
|\det(a_1, a_2,  \ldots , a_k)|$,
and the Lemma is proved.  $\qs$

Returning to the proof of the theorem,
apply the part 1) of the Lemma to
deduce that 
the ray $l(v)$ is contained in one of the cones 
$C_{\FF_i}$, therefore 
one of the 
families $\FF_i$ satisfies the condition 
2) of Theorem \ref{t:exist-regular-bis}.
This family is therefore admissible, 
and $|\det (\FF_i)|<|\det(\FF_0)|$, which contradicts to 
the assumption
that  $|\det (\FF_0)|$ was minimal
among admissible families.  
The theorem is proved. $\qs$

%\newpage

\section{Conical completions and rings of special power series }
\lb{s:twisted}

We will show in this section that
the conical completions of group rings 
in the case of $\xi$-regular families
admit a description in terms of power and Laurent  series
of special type. Let us start
with the simplest case
of free abelian group,
where this description is 
easy to formulate.
For a multi-index $I=(i_1, 
\ldots , i_k)\in \nn^k$
we denote by $t^I$ the monomial 
$t_1^{i_1} \cdot \ldots \cdot t_k^{i_k}$
in variables $t_1, \ldots , t_k$.

\bede\lb{d:special-commut}
We say that a sequence of multi-indices $I_n=(i_1^{(n)}, 
\ldots , i_k^{(n)})\in\zz^k$
{\it strongly converges to $\infty$} 
and we write $I_n \twoheadrightarrow \infty$,
if for every  $j$ with $1\leq j \leq k$ 
we have
$i_j^{(n)} \to \infty$ as $n\to\infty$.
Let $R$ be a commutative ring.
A
series of the form
\bq\lb{f:strongly-conv}
\l=\sum_{n=0}^\infty a_n t^{I_n},
\ \ \ 
a_n \in R, 
\ \ \ 
I_n\in \nn^k,
\ \ 
{\rm and } \ \ 
I_n\tra \infty
\end{equation}
will be called {\it special power series}.
The set of all special power series
will be denoted by $R   \ltrbr t_1, \ldots , t_k\rtrbr$.
It has a natural structure of a ring
and the  inclusions
$R[t_1, \ldots , t_k]\sbs 
R   \ltrbr t_1, \ldots , t_k\rtrbr
\sbs R[[t_1, \ldots , t_k]]
$
are  ring \ho s. $\qt$
\end{defi}
\noi
\bere\lb{r:ser-pol}
Let 
$i\in 
\ldbrack 1, k\rdbrack$.
Any special power series $\l$ can be considered as a power series in 
$t_i$ with coefficients in the ring
$R[t_1, \ldots , t_{i-1}, t_{i+1}, \ldots, t_k]$ 
of polynomials in all other variables. 
That is, we have a natural inclusion
$$ R \ltrbr t_1, \ldots , t_k\rtrbr
\sbs 
\Big(R[t_1, \ldots , t_{i-1}, t_{i+1}, \ldots, t_k]\Big)
[[t_i]].  \hspace{3cm}\qt 
$$
\enre
\noi
Similarly one defines the 
ring of {\it special Laurent series }
with coefficients in $R$. 
The next proposition is obvious;
it will be generalized to non-abelian case later on.

\bepr\lb{p:spec-ser-abel}
Let $\xi:\zz^k\to\rr$
be a monomorphism,
and $\G=(\xi_1, \ldots, \xi_k)$
be a $\xi$-regular family.
Let 
$t_1, \ldots , t_k$
be the basis of $\zz^k$ 
dual to the basis
$(-\xi_1, \ldots, -\xi_k)$
of $\Hom(\zz^k, \zz)$.
Put $\L=\zz[\zz^k]$.
The ring $\Lgam$ 
is isomorphic to the ring of special
power series in the variables $t_1, \ldots , t_k$
with coefficients in $\zz$
and the ring $\Lga$ is
isomorphic to the ring of special
Laurent series in the variables $t_1, \ldots , t_k$
with coefficients in $\zz$. $\qs$
\enpr

\subsection{Twisted special power series}
\lb{su:twi}

Now let us proceed to the non-commutative case. 
Apparently a  first
 example of a   polynomial ring in one variable that 
 does not commute with the elements of the coefficient ring
 ({\it a twisted polynomial ring}) 
 was considered by O. Ore \cite{Ore}, see also  a book of P.M. Cohn
 \cite{Cohn}, \S 2.
 These rings and their 
 generalizations (twisted Laurent 
 extensions of rings, skew power series etc.)
 were thoroughly studied from the point
 of view of their intrinsic structure (see for example 
 \cite{GoldieMichler}), 
 as well as from the point of view of their K-theoretic invariants
 (see for example \cite{FarrellHsiang}).
  However these 
 generalizations still do not cover the algebraic structures
 arising in our present work. So we begin by 
 a brief account of the basic notions of the 
 theory of twisted polynomial rings in the form 
 suitable for our needs.

 Let $A$ be a ring with a unit, $R$ a subring of $A$,
 and $\t_1, \ldots , \t_k\in A$.
 For 
 $I=(i_1, \ldots, i_k)\in \nn^k$
 we denote by 
 $\t^I$ the element $\t_1^{i_1}\cdot\ldots \cdot\t_k^{i_k}$.
 If $\t_i$ are invertible, a similar notation
 will be  used also for   $I\in \zz^k$.
 For $I=(1, \ldots, 1)$ we abbreviate $\t^I$ to $\t$.
 Let $\Sigma= \{\s_1, \ldots, \s_k\}$ 
 be a family of automorphisms of $R$.
 
 \bede\lb{d:tw-poly}
 We say that $A$ is a {\it $\Sigma$-twisted polynomial ring 
 in variables $\t_1, \ldots , \t_k$ with coefficients in $R$}
 (or simply {\it twisted multivariable polynomial ring})
 if
 \been\item
 For every $i,j$ we have
 $\t_i\t_j=r_{ij}  \t_j\t_i$ for some $r_{ij}\in R$.
 \item  For every $i$ and every $r\in R$  we have
 $\t_i r =\s_i(r)\t_i$.
 \item
 $A$ is a free left $R$-module with basis 
 $(\t^I)_{I\in\nn^k}$. 
 \enen 
 These conditions determine the ring $A$ up to an \iso. 
 This ring will be also denoted by $R_\Sigma[\t_1, \ldots , \t_k]$. $\qt$
 \end{defi}
 Observe that the conditions 1) and 2) imply that 
 $R\t^I\cdot R\t^I =R\t^{I+J}$ for all $I, J\in\nn^k$.

 \bere\lb{r:ore}
 Ore's twisted polynomial ring satisfies the conditions
 of the above definition with  $k=1$.
 \enre
 
 The next proposition
  is obvious.
  \bepr\lb{p:ideal}
  Let $A$
  be is a $\Sigma$-twisted polynomial ring 
 in variables $\t_1, \ldots , \t_k$.
 Let $n\in\nn$.
 Then the set $\t^nA$ is a two-sided ideal of $A$
 and the quotient $A/\t^nA$ is a free left $R$-module.
 Its basis is formed by the elements $\t^I$ with $I\in S_n$
 where 
 \bq\lb{f:def-sn}
 S_n=\{ I=(i_1, \ldots, i_k)\in \nn^k~|~
 {\rm \ \ at \ least \ one \ of \ } i_j 
  {\rm \ \ is } \leq n-1 \}.\ \  \ \ \ \qs
  \end{equation}  
\enpr

We will abbreviate $A/\t^nA$ to $A/\t^n$.

\bede\lb{d:tw-Laur}
 We say that $A$ is a {\it $\Sigma$-twisted Laurent  polynomial ring 
 in variables $\t_1, \ldots , \t_k$ with coefficients in $R$}
 (or simply {\it twisted multivariable L-polynomial ring})
 if
 \been\item
 Elements $\t_i$ are invertible in $A$.
 \item 
 For every $i,j$
 the commutator $\t_i\t_j\t_i^{-1}\t_j^{-1}$
 is in $R$.
  \item  For every $i$ we have
 $\t_i r \t_i^{-1}=\s_i(r)$.
 \item
 $A$ is a free left $R$-module with basis 
 $(\t^I)_{I\in\zz^k}$. 
 \enen 
 
 \noindent
  These conditions determine the ring $A$ up to an \iso. 
 This ring will be also denoted by $R_\Sigma[\t^\pm_1, \ldots , \t^\pm_k]$. $\qt$
  \end{defi}
   \noindent
Similarly to the case of twisted polynomial rings  we have 
 $R\t^I\cdot R\t^I =R\t^{I+J}$ for all $I, J\in\zz^k$.
%Observe that the properties 2) and 3) of this definition
%are equivalent to the properties 1) resp. 2) of the 
%definition \ref{d:tw-poly}

\begin{exam}\lb{e:grr-l}
A basic example of a twisted Laurent polynomial ring
 arises as follows. Let $G$ be a group  endowed with 
 an epimorphism $\z:G\to \zz^k$, let $H=\Ker\z$.
 Choose some free generators 
 $t_1, \ldots, t_k$
 of $\zz^k$
  and choose 
 any elements $\t_i\in G$ \sut~ $\z(\t_i)=t_i$.
 Denote by $\s_i$ the automorphism of $\zz H$
 defined by $\z_i(x)=\t_ix\t_i^{-1}$, and put $\Sigma=(\s_1, \ldots, \s_k)$.
 Then the natural \ho~ 
 $\zz H_\Sigma [\t^\pm_1, \ldots , \t^\pm_k]\to \zz G$ 
 is an \iso, so that the 
 group ring 
 $\zz[G]$ is isomorphic to a twisted Laurent polynomial ring
 in $k$ variables.
 \end{exam}

The next proposition   is obvious.
  \bepr\lb{p:poly-in-Laur}
  Let $A$ be  a $\Sigma$-twisted L-polynomial ring 
 in variables $\t_1, \ldots , \t_k$. 
 Denote by $A_0$ the left $R$-submodule
  of $A$ generated by $\t^I$ with $ I\in\nn^k$.
  Then $A_0$ is a subring of $A$ and it is a 
  $\Sigma$-twisted polynomial ring in variables $\t_1, \ldots , \t_k$. $\qs$
    \enpr
    
    Let $A$ be a $\Sigma$-twisted Laurent polynomial ring 
 in variables $\t_1, \ldots , \t_k$. 
 For $j\in \ldbrack 1, k\rdbrack$
 denote by $A_j$ the free left 
 $R$-submodule generated by the elements $\t^I$ where 
 the multi-index $I=(i_1, \ldots, i_k)\in\zz^k$
 satisfies the 
 condition $i_j\geq 0$. It is clear that $A_j$ is a ring. 
 Let $n\in\nn$. 
  The inclusion $A_0\hookrightarrow A_j$
induces the map of quotient 
 rings 
 $$
 J_{n,j}:A_0/\t^n \to A_j/\t_j^n, 
 $$
 Denote by $J_n$ the direct sum of these maps
 $$
 J_n 
 =
 \bigoplus_j J_{n,j} : A_0/\t^n \to \bigoplus_j \Big( A_j /\t_j^n\Big),
 \ \  j\in \ldbrack 1, k\rdbrack. 
 $$
 
 \bepr\lb{p:inject}
 For any $n$ the map 
 $J_n$ is injective.
 \enpr
 \Prf
 Let $x\in A_0/\t^n$, 
 write $x=\sum_I x_I\t^I$; 
 here the sum ranges over $I\in S_n$.
  Assume that $x_I\not=0$ for some $I= (i_1, \ldots, i_k) $.
 By definition of $S_n$ there is an integer 
 $r\in \ldbrack 1, k \rdbrack$
 \sut~ $i_r<n$.
 Then the image of $x$ in $A_{r}/\t^n_{r}$ is non-zero. $\qs$
 
 Let us proceed to completions of
 twisted polynomial rings.
 
 \bede\lb{d:inv-lim-twi}
 Let $A$ be a twisted  polynomial ring 
 in variables $\t_1, \ldots , \t_k$. 
 Consider  the sequence 
 $$
 A \leftarrow A/\t \leftarrow A/\t^2 \leftarrow \ldots 
 $$
of natural projections and denote by $\wh A$
 its inverse limit.
  For $j\in \ldbrack1, k \rdbrack$ 
  consider  the sequence 
  $$
 A_j \leftarrow A_j/\t_j \leftarrow A_j /\t_j^2 \leftarrow \ldots 
 $$
 of natural projections
 and denote by 
  $\wh {A_j}$
 its inverse limit.
  We have then a natural ring \ho~ 
 $$
 \wh J_{j} : \wh A \rightarrow \wh {A_j}.\ \ \ \ \ \ \ \ \ \ \qt
 $$
  \end{defi}
 
 \noindent
 The next proposition is obvious.
 \bepr\lb{p:inj-compl}
 Let $A$ be a twisted  polynomial ring. Then 
the \ho~ $
 \wh J_{j}$ 
 is injective 
for every $j$.
$\qs$
\enpr
 
The completion of a twisted polynomial ring
can be seen as a ring of special power series 
(in non-commuting variables).

\bede\lb{d:convv-ind}
Let $A$ be a $\Sigma$-twisted  polynomial ring 
 in variables $\t_1, \ldots , \t_k$. 
%We say that a sequence of multi-indices $I_n=(i_1^{(n)}, 
%\ldots , i_k^{(n)})\in\zz^k$
%{\it strongly converges to $\infty$} 
%and we write $I_n \twoheadrightarrow \infty$,
%if for every  $j$ with $1\leq j \leq k$ 
%we have
%$i_j^{(n)} \to \infty$. 
Consider the set of all 
series of the form
\bq\lb{f:strongly-convv}
\l=\sum_{n=0}^\infty a_n\t^{I_n},
\ \ \ 
a_n \in R, 
\ \ \ 
I_n\in \nn^k,
\ \ 
{\rm and } \ \ 
I_n\tra \infty.
\end{equation}

This set has a natural ring structure
determined by properties 1) and 2) of Definition \ref{d:tw-poly}.
We call it 
{\it $\Sigma$-twisted special power series ring in variables
$\t_1, \ldots , \t_k$} or simply
{\it twisted sp-series ring}
(in order to distinguish it from the usual 
power series ring).
We denote it by 
$R_\Sigma   \ltrbr \t_1, \ldots , \t_k\rtrbr$;
the inclusion
$R_\Sigma[\t_1, \ldots , \t_k]\sbs R_\Sigma   \ltrbr \t_1, \ldots , \t_k\rtrbr$
is a ring \ho.
$\qt$
\end{defi}
\noindent
The next proposition is obvious.

\bepr\lb{p:seriess}
Let $A$ be a $\Sigma$-twisted  polynomial ring 
 in variables $\t_1, \ldots , \t_k$. 
 There is a natural isomorphism
 \begin{flalign*}
 &&
\wh A
 \approx 
 R_\Sigma   \ltrbr \t_1, \ldots , \t_k\rtrbr.  && \qs
 \end{flalign*}
 \enpr

 %Similarly we define the special Laurent
 %power series ring.
 \bede\lb{d:laur-series}
 Replacing $\nn^k$ by $\zz^k$ in the formula
 \rrf{f:strongly-convv}
 we obtain the definition
 of 
 {\it $\Sigma$-twisted special Laurent series ring in variables
$\t_1, \ldots , \t_k$} or simply
{\it twisted sl-series ring}.
We denote it by 
$R_\Sigma   ((( \t_1, \ldots , \t_k)))$;
the inclusion
$
R_\Sigma   \ltrbr \t_1, \ldots , \t_k\rtrbr
\sbs
R_\Sigma   ((( \t_1, \ldots , \t_k)))
$
 is a ring \ho.
 \end{defi}

\subsection{Conical completions of group rings  }
\lb{su:cc-xi-reg}

Now we can give an interpretation
of conical completions in terms of special  power series.
%Let us begin by considering group rings.

Let $G$ be  group. We denote by $\L$ the group 
ring $\zz G$. Let $\z:G\to \zz^k$ be an epimorphism,
and
$\xi : G\to\rr,\ \  \xi_1, \ldots, \xi_k : G\to\zz$
\ho s factoring through $\z$,
that is, 
there are \ho s $\bar \xi:\zz^k\to\rr, \      \bar \xi_i:\zz^k\to\zz$
\sut~ $\xi_i=\bar\xi_i\circ \z$, and $\xi=\bar\xi\circ\z$. 
\bede\lb{d:reg_g}
The family $\G=\{\xi_1, \ldots, \xi_k\}$ 
is called {\it $\xi$-regular}
if the family 
$
\{\bar\xi_1, \ldots, \bar\xi_k\}$
is $\bar\xi$-regular. $\qt$
\end{defi}
Denote by 
 $(t_1, \ldots, t_k)$
the basis in $\zz^k$ dual to
the basis $(-\bar\xi_1, \ldots, -\bar\xi_k)$ of $\Hom(\zz^k, \zz)$.
Pick any elements $\t_1, \ldots,\t_k\in G$
with $\z(\t_i)=t_i$.
Let 
$$H_i=\Ker\xi_i{\rm\ \  and }
\ \ H=\bigcap_i H_i = \Ker (\z~:~ G\to\zz^k).$$
The subgroups $H_i$ and $H$ are normal in $G$,
and the element $\t_i$ determines an automorphism $\s_i$ of $H_i$ and $H$ 
as follows: $\s_i(x)= \t_i x \t_i^{-1}$.
Let $\t= \t_1\cdot \ldots \cdot \t_k$.
\bede\lb{d:lambdas-and-gammas}
Let $n\in\nn$.

1) \ \ Put $$\lin=\left\{\l\in\L~|~ \supp\l\sbs \xi_i^{-1}\Big(]-\infty, -n]\Big)\right\}.$$
 We abbreviate $\L_{i,0}^-$ to $\linol$.
Then $\linol$ is a ring with a unit, 
and  $\lin$ is a two-sided  ideal of
$\linol$, generated as a left (or right) ideal by $\theta_i^n$.
The quotient ring $\linol\big/\lin$
will be denoted by $\linn$.

2) Put $\lign = \bigcap_i \lin$.
We abbreviate $\L_{\G,0}^-$ to $\lgamin$.
Then $\lgamin$ is a ring with a unit, and 
$\lign$ is a two-sided principal ideal of
$\lgamin$.
%, generated by $\theta^n$.
The quotient ring 
$\lgamin\big/\lign$
will be denoted by $\lgann$.

3) The inclusion 
$\lgamin \hookrightarrow \linol$
induces \ho s
$$
\wi J_{n,i} : \lgann \to \linn.
$$
We denote by 
$\wi J_{n}$
the direct sum of these maps
 \begin{flalign*}
 && \wi J_n 
 =
 \bigoplus_i \wi J_{n,i} : \lgann \to \bigoplus_i \linn,
 \ \  i\in \ldbrack 1, k\rdbrack. && \qt
  \end{flalign*}

\end{defi}
The natural isomorphism from the example 
\ref{e:grr-l}
induces the isomorphisms of the next proposition.
\bepr\lb{p:polynoms-lambdas}
Let $n$ be a natural number and $i\in \Lbrack 1, k \Rbrack$.
Put $\Sigma=\{\s_1, \ldots , \s_k\}$.
We have natural isomorphisms
\begin{equation}\lb{f:iso1}
\linol \approx (\zz H_i)_{\s_i}[\t_i], \ \  \linn\approx   
\big((\zz H_i)_{\s_i}[\t_i]\big)/\t_i^n;  
\end{equation} 
\begin{equation}\lb{f:iso2}
%\lin\approx (\zz H_i)[\t^n_i];  \\
\lgamin \approx (\zz H)_\Sigma[\t_1, \ldots, \t_k];  
\end{equation}
\begin{equation}\lb{f:iso3}
\lgann \approx \Big((\zz H)_\Sigma[\t_1, \ldots, \t_k]\Big)\Big/\t^n. \ \ \ \ \ \qs
\end{equation}
\enpr 

The next Corollary follows from \ref{p:inject}. 

\beco\lb{c:monomorphisms}
The direct sum $\wi J_n$ of the maps $
\wi J_{n,i} : \lgann \to \linn $
 \begin{flalign*}
 && \wi J_n 
 =
 \bigoplus_i \wi J_{n,i} : \lgann \to \bigoplus_i \linn,
 \ \  i\in \ldbrack 1, k\rdbrack. && 
  \end{flalign*}
is injective. $\qs$
\enco

Now we can proceed to conical completions.
We will use the following abbreviations:
$$
\Lii= \Lxii, 
\ \ \ \ 
\Liimin= \Lxiimin.
$$

The next proposition is obvious.
\bepr\lb{p:identif}
The isomorphism 
\rrf{f:iso2}
extends to  natural isomorphisms
$$\Lgam \approx (\zz H)_\Sigma
\ltrbr \t_1, \ldots , \t_k\rtrbr,
\ \ 
\Lga \approx (\zz H)_\Sigma
((( \t_1, \ldots , \t_k ))).
$$

The isomorphisms
\rrf{f:iso1}
induce natural isomorphisms
$$
\Lii \approx (\zz H_i)_{\s_i}((\t_i)), \ \  
\Liimin \approx (\zz H_i)_{\s_i}[[\t_i]], $$
$$
\Liimin/\t^n\Liimin \approx\linn\approx   
\big((\zz H_i)_{\s_i}[\t_i]\big)/\t_i^n;  $$
\enpr
\noi
Observe that the inclusion
$$
\Lgam \sbs \Liimin 
$$
is 
an analog of the injective map 
$\wh J_{j}$
from Proposition \ref{p:inj-compl}.

\bele\lb{l:injec}
The natural \ho s
$$
\l^-_i: \Lgam \to \Liimin, 
\ \ \ \ 
\l_i: \Lga \to \Lii
$$
are injective for every $i$. $\qs$ 
\enle

%%%%%%%%%%%%%%%%%

%\newpage

\section{Conical refinement of the Novikov complex: the statement of the main theorem}
\lb{s:con_Nov}

Let 
$\o$ be a Morse form on a closed connected \ma~ $M$.
Denote by $\xi\in H^1(M,\rr)$ its de Rham cohomology
class. 
Let $p:\wh M\to M$ be a regular covering 
with a structure group $G$, such that 
$[p^*(\o)]=0$.
Then the cohomology class $\xi$ can  be considered as a \ho~ $\xi : G\to\rr$.
Let us factor it as follows:

$$\xymatrix{
G \ar[dr]_\zeta\ar[rr]^{\xi} & & \rr \\
& \zz^k \ar[ur]_{\bar\xi} & 
}
$$
where $\bar\xi$ is a \mo~ and $\z$ is an \epi.
Recall that $k$ is called the {\it irrationality degree}
of $\o$. 
Let $\L=\zz G$.
Denote by $Z(\o)$ the set of all zeros of $\o$ and by $Z_r(\o)$ 
the set of all zeros of $\o$ of index $r$.

\beth\lb{t:main_th}
Let $v$ be a transverse
$\o$-gradient. 
Then there is a $\xi$-regular family 
$\G=\{\xi_1,\ldots, \xi_k\}$
of \ho s $G\to\rr$
and a  chain complex 
$\NN_*(\o,v, \Gamma)$ 
of free $\Lga$-modules,
freely generated in degree $r$ by 
$Z_r(\o)$,  and chain equivalent 
to  $\SS_*(\wh M)\tens{\L}\Lga   $.
\enth
The proof is done in  Section 
\ref{s:prf_main}.
%\ref{s:prf_general}.
For the case $k=1$ the theorem was proved in \cite{P-o}.
The next section contains a brief outline of this proof.
The geometric part of the proof of Theorem 
\ref{t:main_th} 
follows the lines 
of \cite{P-o}.

\section{The case of circle-valued Morse functions}
\lb{s:circle}

Assume that $k=1$, that is, $\o=\l\cdot df$, where 
$\l>0$ and 
$f:M\to S^1$ is a Morse function non-homotopic to zero.
Then $\xi=\l [f]$ and $[f]\in H^1(M,\zz)$.
We can assume that $\l=1$ and that $\xi$ is indivisible.
%This particular case of the main theorem was done in the author's
%work \cite{P-o}. The present section contains 
%a brief outline of this paper. 
In this section we give only brief 
indications for the proofs of Propositions and Lemmas, 
referring to \cite{P-o} for full proofs.

Consider the infinite cyclic covering
$p: \ove M \to M$ induced from the universal covering
$\rr\to S^1$ by the map $f:M\to S^1$.
Lift the function $f:M\to S^1$ to a Morse function 
$\bar f :\ove M \to \rr$. We have  a commutative diagram 

$$\xymatrix{
\wh M \ar[rr]^Q \ar[dr] &  & \ove M \ar[dl]^p\\
& M & 
}
$$
\noi
Put $\wh f=\bar f\circ Q$. Assume that $0$ is a regular value
of $\bar f$ and let 
\bq\lb{f:set-Y}Y=\{x~|~ \wh f (x) \leq 0\},
\ \ \ \ 
Y_n=\{x~|~ \wh f (x) \leq -n\}.
\end{equation}
Denote by $t$ the generator of the structure group $\approx \zz$ 
of the covering $\ove M \to M$, satisfying $\bar f(tx)=\bar f(x) -1$. In this section
we will abbreviate 
$\Lxim$ to 
$\wh\L^-$ and 
$\Lxi$ to 
$\wh\L$.
Let 
\begin{gather*}
\wh\L^-_n 
=
\{\l\in \Lxim ~|~ \supp\l\sbs \xi^{-1}(]-\infty, -n])\},\\
 \L^{(n)} = \wh\L^-/\wh\L^-_n.
\end{gather*}

\subsection{Definition of the Novikov complex}
\lb{su:def-plan}

We use the classical Witten's reformulation of the Morse theory
\cite{Witten1982}.
Let $v$ be a transverse $f$-gradient. Its lifts to $\wh M$ and $\ove M$ 
will be denoted by the same letter $v$.
For $a\in Crit_k(\wh f), \ b\in Crit_{k-1}(\wh f)$,
denote  by $\G(a,b)$ the set of flow lines of $v$ in $\wh M$ 
joining $a$ with $b$. Applying usual Morse theory to the cobordism
$\bar f^{-1}([A, B])$, where $\wh f(a), \wh f(b)\in ]A, B[$, 
it is easy to see that $\G(a,b)$
is finite. For each $p\in Crit(f)$ choose an orientation 
of the stable manifold of $p$ \wrt~ $v$,
this choice induces orientations for stable manifolds of all
critical points of $\wh f:\wh M\to \rr$. 
Then each flow line $\g\in \G(a,b)$
is endowed with a sign 
$\ve(\g, v)\in \{-1, 1\}$.
Let $$n(a,b)= \sum_{ \g\in \G(a,b)} \ve(\g).$$
For every critical point $p\in Crit(f)$ let $\bar p$ be the lift 
of $p$ to the cobordism $\bar f^{-1}([-1,0])$.
Put
$$
N(p,q)=\sum_g n(\bar p, \bar q g)\cdot g \in \lLL.
$$

\noindent
It is clear that 
$N(p,q) \in \lilim$.
Let $\NN_k^-$ be the free $\lilim$-module
generated by $Crit_k(f)$. Put
$$
\pr_k(p)=\sum_{q\in Crit_{k-1}(f)} N(p,q)q;
$$
we obtain a \ho~ 
$\pr_k:\NN_k^-\to \NN_{k-1}^-$.
It is not difficult to prove that $\pr_{k-1}\circ\pr_k =0$. Indeed,
it suffices to check that the image of 
$\pr_{k-1}\circ\pr_k (p)$ 
in $\lilin$ 
vanishes for every natural number $n$ 
and every critical point $p$ in $Crit_k(f)$,
and this is proved by applying the classical Morse theory
to the cobordism 
$
f^{-1}([-n, 0])$.

We obtain therefore a chain complex $\nin$ 
of free $\lilim$-modules, generated in degree $k$ by $Crit_k(f)$.
Let $\t\in G$ be any element in $G$, \sut~ $\z(\t)=t$.
To relate the chain complex $\nin$ to the completion of $\SS_*(\wh M)$
we first 
construct a chain equivalence $J_n$ between
$\SS_*(Y, Y_n)$ and  
$\nin \big/ \t^n\nin\approx \NN^-_* \tens{\lilim} \lilin$
(see Subsection \ref{su:t-ord}).
Moreover  these chain equivalences can be chosen to be compatible 
 with each other for different values of $n$.
In the second part of the proof
(see Subsection \ref{su:limit})
we use  this compatibility 
to pass to the limit as $n\to  \infty$,
and this determines the required chain equivalence.

\subsection{$t$-ordered Morse functions}
\lb{su:t-ord}

In this subsection we outline the construction of the chain equivalence 
$J_n$. Let $W_n=\bar f^{-1}(]-n, 0])$, and 
$\wh W_n=Q^{-1}(W_n)$. Then $W_n$ is a cobordism,
its boundary is the disjoint union of two manifolds
$\pr_0 W_n=\bar f^{-1}(-n)$, and 
$\pr_1 W_n=\bar f^{-1}(0)$.
Observe  that the chain equivalence of 
$\D_*(Y)/\t^n$ and $\NN_*^-/\t^n$ as chain complexes over $\zz$ follows from the classical 
Morse theory, since $\NN_*^-/\t^n$ is the Morse complex of the Morse function
$\wh f:\wh W_n\to\rr$ and $\D_*(Y)/\t^n$ is 
isomorphic 
to 
$\D_*(\wh W_n, \pr_0\wh W_n)$.
To construct a chain equivalence
respecting the $\lilin$-structure, we introduced in \cite{P-o}, Lemma 5.1
the notion of $t$-ordered Morse function.

\bede\lb{d:t-ordered}
Let $m=\dim M$.
A Morse function $\Psi:W_n\to [-1,m]$
on the cobordism $W_n$ is called {\it $t$-ordered Morse function adjusted to $f$ and $v$},
(or just {\it $t$-ordered} if \noconf) if

1) $v$ is a $\Psi$-gradient,

2) $\Psi(tx)< \Psi(x)$ for every $x$,

3) For every critical point $x$ of $\Psi$ of index $r$ we have
$r-1< \Psi(x) < r$.
\end{defi}

The existence of $t$-ordered functions was proved in \cite{P-o}.
\footnote{There is a minor difference between
Definition \ref{d:t-ordered}
and the definition from \cite{P-o},
namely, the image of $t$-ordered function in \cite{P-o} 
can be any closed interval of $\rr$.}
Pick a $t$-ordered Morse function $\Psi$ and extend it 
to the whole of $\ove M$
in such a way that 
$\Psi^{-1}([-1, m]) =W_n$.
Put $\wh \Psi= \Psi\circ Q$ and let 
$Z_r=\{x~|~ \wh \Psi(x)\leq r\}$,
so that we have
\bq\lb{f:filt-Y}
Y_{n}= Z_{-1}\sbs Z_{0}\sbs \ldots Z_{m} = Y.
\end{equation}
Then the pair $(Y,Y_n)$ is filtered by pairs
$(Z_r, Y_n)$, where $r$ ranges over integers in $ [-1, m]$. 
The homology of the pair 
$(Z_r,Z_{r-1})$
is computed via the classical Morse-theoretic
procedure. Namely, for every critical point $p$ of $f$ 
in $Z_r\sm Z_{r-1}$
the intersection $D_p$ of the stable manifold of $p$ 
with $Z_r\sm \Int Z_{r-1}$
is homeomorphic to the $r$-dimensional closed disc 
and determines therefore a homology class
$\D_p\in H_r(Z_r,Z_{r-1})$.
Then $H_r(Z_r,Z_{r-1})$
is a free abelian group generated by the classes $\D_p$
where $p$ ranges over the set of critical points of $\wh\Psi$
belonging to $Z_r\sm Z_{r-1}$.
Denote this group by
$\FF_r^{(n)}$.
%=H_r(Z_r,Z_{r-1})$.
The boundary operator of the exact sequence of the triple
$(Z_r,Z_{r-1},Z_{r-2} )$
endows the graded module $\FF_*^{(n)}$
with the structure of a chain complex.

\bele\lb{l:isoo}
The map $p\mapsto \D_p$
induces an isomorphism
$$
J^{(n)}: \NN_*^-/\t^n \to \FF_*^{(n)}
$$
of chain complexes 
over $\lilin$.
Therefore $H_r(Z_r,Z_{r-1})$
is a free $\lilin$-module with basis
$Crit_r(f)$.
%$(\D_p)_{p\in Crit_r(f)}$.
\enle
The proof of the Lemma is based on the equality
$\D_{pg}=\D_p\cdot g$
which holds for every $g\in\lilin$. This equality 
follows from the fact that $\Psi$ is a $t$-ordered 
Morse function. $\qs$

Thus the filtration of the $\lilin$-chain complex $\SS_*(Y, Y_n)$
by subcomplexes $\SS_*(Z_r, Y_n)$ is {\it cellular} 
(or {\it good} in the terminology of \cite{P-o}).
Applying the classical method of computing the homology
of CW complexes from the complex of cellular chains
we deduce that the homology of $\FF_*^{(n)}$
is isomorphic to $H_*(Y, Y_n)$.
Moreover one can prove that there is a canonical chain equivalence
between $\FF_*^{(n)}$ and $S_*(Y, Y_n)$
(see \cite{P-o}, \S 3).
We obtain therefore a chain equivalence of $\lilin$-complexes
$$
J_n\ : \ \NN_*^-/\t^n \arrr \sim \SS_*(Y, Y_n).
$$
\bele\lb{l:square-commut}
The following diagram is homotopy commutative
(the horizontal arrows below are natural projections).
$$
\xymatrix{ 
\NN_*^-/\t^n  \ar[d]^{J_n}    &  \NN_*^-/\t^{n+1} \ar[l]\ar[d]^{J_{n+1}}\\
\SS_*(Y, Y_n)  & \SS_*(Y, Y_{n+1})\ar[l]
}
$$
\enle
The basic observation for the proof of this lemma
is that one can choose the ordered Morse function
on the cobordism $W_{n+1}$ to be compatible with 
a given ordered Morse function on $W_n$ (see \cite{P-o}),
so that the projection map $(Y, Y_{n+1}) \to (Y, Y_{n})$
preserves filtrations. $\qs$

%%%%%%%%%%%%%%%%%%%%%%%%%%%%

\subsection{Inverse limit of the chain equivalences}
\lb{su:limit}

Replacing the singular chain complex 
$\SS_*(Y, Y_n)$
by simplicial chain complex
$\D_*(Y, Y_n)$
we obtain a chain equivalence 
$\nin \big/ \t^n\nin
\arrr {I_n}
\D_*(Y)/\t^n\D_*(Y)$, such that
%%Moreover the chain equivalence can be chosen is such a way
%that for every $n$ 
the following diagram is homotopy commutative.

\begin{equation}\lb{f:diagr}
\xymatrix{ 
\nin/\t^n \ar[d]^{I_{n}}& \nin/\t^{n+1}\ar[l] \ar[d]^{I_{n+1}} \\
\D_*(Y)/\t^n &   \D_*(Y)/\t^{n+1}  \ar[l]
}
\end{equation}
%This part of the proof is based on the theory 
%of $t$-ordered Morse functions, introduced 
%and developed in \cite{P-o}.
Recall  that 
$
\wh\D_*(Y)\approx \liminv\ \wh\D_*(Y)/\t^n,
$ and 
$
\nin \approx \liminv\ \nin/\t^n.
$
Using the diagram \rrf{f:diagr}
we construct a chain complex 
$Z_*$ and homology equivalences
$$
\wh\D_*(Y)\arr Z_* \arl \nin$$
of $\lilim$-complexes
and deduce from it 
a chain equivalence 
\bq\lb{f:equiv-}
\wh\D_*(Y) \sim \nin.
\end{equation}
This last argument 
(see \cite{P-o}, \S 3, part B)
is purely algebraic.
Observe that the diagram 
\rrf{f:diagr} 
does not imply immediately the required
chain equivalence, since it is commutative only up to homotopy.

Take the tensor product 
of the chain equivalence \rrf{f:equiv-}
by $\lili$  over $\lilim$,
replace the simplicial chain complex on the left 
by the singular chain complex, and  
we obtain the chain equivalence sought.

%\newpage

\section{Proof of the main theorem }
%for the case $G= \zz^k$ and $k={\rm irr}\ \o$}
\label{s:prf_main}

%In this particular case we have $\wh M = \wi M$, and
%the structure group of the covering $p:\wi M \to M$ 
%is isomorphic to
%$\zz^k$. 

Let $\O=(\o_1, \ldots , \o_k)$ be a $(\o,v)$-regular
family; put $\xi_i=[\o_i]$
and let $\G=(\xi_1, \ldots , \xi_k)$.
We will construct
the Novikov complex $\NN_*(\o,v, \G)$ defined over $\Lga$.
For every $i$ we have a commutative diagram of 
coverings
$$\xymatrix{ \wh M \ar[ddr]_{p}\ar[r]^{Q} &
\wi M \ar[dd]^{\wi p} \ar[dr]^{Q_i} &    \\
&   &\quad \ove M_i  \ar[dl]^{p_i}   \quad\label{f:determ} \\
& M             & 
}
$$

\noindent
Here $p_i:\ove M_i\to M$ 
is the infinite cyclic covering corresponding to
the integer cohomology class $\xi_i\in H^1(M,\zz)$, the map 
 $\wi p : \wi M 
\to M$ is a covering with structure group $\zz^k$,
and $Q_i$ is 
%The covering $p$ factors through
a regular covering  with structure group $\zz^{k-1}$.
%The ring $\wh\L_{\xi_i}$ in this case is easy to describe.
Let  $(t_1, \ldots , t_k)$ be the  $\zz$-basis of $\zz^k$,
dual to the  $\zz$-basis $(-\xi_1, \ldots , -\xi_k)$ of $\Hom(\zz^k, \zz)$.
Choose $\t_i\in G$ \sut~ $\z(\t_i)=t_i$.
Then 
$\wh\L_{\xi_i}$
is isomorphic to the ring $R_i((\t_i))$ of
twisted
Laurent series in $\t_i$ with coefficients in
the group ring $R_i$ of the group $\Ker\xi_i$.
The ring $\Lga$ 
is isomorphic to
$R(((\t_1, \ldots, \t_k)))$,
where 
$R=\zz[\Ker \xi]$
and the ring
$\Lgam$ 
is isomorphic to
$R\ltrbr \t_1, \ldots, \t_k\rtrbr$.
%Then  
%$\wh\L_{\xi_i}$
%is isomorphic to the ring $R_i((t_i))$ of
%Laurent series in $t_i$ with coefficients in 
%the Laurent polynomial ring $R_i$
%in $(k-1)$ variables 
%$t_1, \ldots , t_{i-1}, t_{i+1}, 
%\ldots, t_k$.
Let $f_i:M\to S^1$ be a circle-valued Morse function
\sut~ $df_i=\o_i$;
then $v$ is also an $f_i$-gradient.

\subsection{Construction of the Novikov complex}
\lb{su:constr_nov}

To construct the Novikov complex 
we follow the same schema as in Section \ref{s:circle}.
For every zero $a$ of $\o$
choose a lift $\wh a $ of $a$ to $\wh M$.
The lift of the vector field $v$ to $\wh M$
will be denoted by the same letter $v$.
Let $a,b\in Z(\o)$ with $\ind a = \ind b +1$. Let $g\in G$.
Denote by 
$\G(a,b;g)$
the set of all flow lines  of $v$
in $\wh M$ joining $\wh a $ with  $\wh b\cdot g $.
Assume that $v$ is a transverse $\o$-gradient.
For each $a\in Z(\o)$ choose an orientation
of the stable manifold 
of the point $a$ with respect to the flow 
generated by $v$.
A standard procedure from Morse-Smale 
theory associates to each $\g\in \G(a,b;g)$ 
a sign $\ve(\g, v)\in \{1, -1\}$.

\bele\label{l:stepA}
The set $\G(a,b;g)$ is finite for every $g\in G$.
\enle
\Prf
It follows immediately from Lemma 2.1
of \cite{P-o}, part (1),
applied to the circle-valued Morse function 
$f_i:M\to S^1$ (where $i$ is
any integer
in $[1,k]$), and the regular covering 
$\wh M\to M$. $\qs$
%having  the property $[p^*(\o_i)]=0$.

\noi
This lemma enables us  to define the Novikov incidence coefficient
$N(a,b; v)$ as follows. Put
$$ n(a,b;g) = \sum_{\g\in \G(a,b;g)}\ve(\g, v)\in\zz, \ \ {\rm and } $$
\bq\lb{d:inc_coeff}
N(a,b;v) = \sum_{g\in G} n(a,b;g)g \in \lL.
\end{equation}

\bele\lb{l:stepB}
We have $N(a,b;v)\in \Lga$.
\enle
\Prf
The  incidence coefficient $N(a,b;v)$ 
belongs  to the ring 
$\wh\L_{\xi_i}$ by Lemma 2.1
of \cite{P-o}, part (2).
This holds for every $i$ therefore we have
$N(a,b;v) \in \Lga= \bigcap_i \wh\L_{\xi_i}$,
and the assertion of the lemma follows.
$\qs$
\pa
Let 
$\NN_r$
be the free $\Lga$-module
generated by $Z_r(\o)$.
Using Lemma \ref{l:stepB} we  define 
a homomorphism  
$\NN_r
\to 
\NN_{r-1}$
by
\begin{equation}\lb{f:d_r}
\pr_r a = \sum_b N(a,b;v) b.
\end{equation}
\bele\label{l:stepC}
We have $\pr_r \circ \pr_{r+1} =0$ for every $r$,
so that the graded 
$\Lga$-module $\NN_*$ endowed with the operator $\pr$
is a chain complex.
\enle
\Prf 
The ring $\Lga$ is a subring of $\wh\L_{\xi_i}$
for every  $i$, 
and the module 
 $\NN_*$
 is a submodule of the Novikov complex
  $\NN_*(\o_i,v)$ of the rational Morse form
  $\o_i$,
so the assertion of the 
lemma follows from Theorem 2.2 part (1) of \cite{P-o}.
$\qs$ 

\subsection{Truncated Novikov complexes}
\lb{su:truncated}

Now
we can begin   the 
%central part of the proof 
%of Theorem \ref{t:main_th}, namely 
%to 
construction of a chain equivalence 
between $\NN_*$ and 
$ \SS_*(\wh M) \tens{\L}  \Lga$. 
%This construction 
%occupies the rest of this section.
Observe that a direct application
of the results of \cite{P-o}
gives only a chain equivalence of chain complexes
over the ring $\wh\L_{\xi_i}$ for each $i$.
We start by constructing chain equivalences between certain truncated versions
of the Novikov complex and the singular chain complex.
%The chain complex $\NN_*$ is chain equivalent
%to $\Lga\tens{\L} C_*(\wh M)$. }
Lift the maps 
$f_i:M\to S^1$ to functions $\bar f_i:\ove M_i \to \rr$ and put
$\wh f_i=\bar f_i\circ Q_i\circ Q:\wh M \to \rr$.
Let 
$$g=\max (\wh f_1, \ldots, \wh f_k);$$
this is a continuous function 
$\wh M \to \rr$.
We can assume that $0$ is a regular value for every function $\wh f_i$,
and that the chosen lifts of zeros of $\o$ to $\wh M$ 
are in the set 
$\{x\in \wh M ~|~ \wh f_i(x)\in ]-1, 0[ \ \ \forall i \}$. 
\ Denote the set of these lifts by $K$ and the 
set of the lifts of zeros of index $r$ by $K_r$. 
Put 
$$
Y=\{x\in \wh M ~|~ g(x)\leq 0\}. \ \ 
$$
Denote by $\NN^-_*$
the abelian subgroup of all formal 
linear combinations $\l$ (infinite in general)
of zeros of $p^*(\o)$ 
belonging to $Y$, and subject to the following condition $(\CC)$:
\pa
$(\CC)$: \ \ For every $C\leq 0$ and every $i$ 
the set 
$\supp \l \cap \wh f_i^{-1}([C, 0])$ is finite.
\pa

\noi
Using the fact that the family $(\o_1, \ldots , \o_k)$ is $(\o,v)$-regular
it is easy to see that 
$\NN_*^-$ has a natural structure of a 
graded free $\Lgam$-module
generated by $K$.
We endow it with the structure of a chain complex over $\Lgam$ using 
the formula \rrf{f:d_r}. We have then 
$$
\NN_*\approx   \NN_*^-\tens{\Lgam} \Lga.
$$
\noindent
Put 
$$
Y_n=\{x\in \wh M ~|~ g(x)\leq -n\},
$$
let $\NN_*^{n}$ be the subcomplex of $\NN^-_*$
formed  by all zeros of $p^*(\o)$ in $Y_n$. Put
$$ \NN_*^{(n)} = \NN^-_*/\NN_*^n
, 
\ \ \ 
{\rm then }
\ \ \ 
\NN^-_* = \liminv\   \NN_*^{(n)}.$$
%Let $\FF=(\wh f_1, \ldots, \wh f_k)$.
%Following the method of \cite{P-o}
%we will do the proof in 3 steps.
\bepr\lb{fl:kvadrat}
 For every $n$ there is  a chain equivalence 
$J_n:\NN_*^{(n)}\arrr \sim \SS_*(Y, Y_n)$
of $\Lgam$-complexes, 
\sut~ the following diagram 
%$(\Delta_n)$ 
is chain homotopy commutative

\bq\lb{f:kvadrat}
\xymatrix{ 
\NN_*^{(n)} \ar[d]^{J_n}    & \NN_*^{(n+1)} \ar[d]^{J_{n+1}} \ar[l] \\
\SS_*(Y, Y_n)  & \SS_*(Y, Y_{n+1})  \ar[l]
}
\end{equation}
(where the horizontal arrows are the natural projections).
\enpr
%This step requires a version of equivariant Morse theory,
%and goes along the lines of \cite{P-o}, with  a minor modification 
%of the arguments therein.
\Prf
For $i\in \ldbrack1, k \rdbrack$
Let $\tau_i$ be the generator of the structure group $\approx \zz$
of the covering $\ove{M}_i \to M$, such that $\ove{f}_i(\tau_i x) = \ove{f}_i(x)-1$.
%For each $i$ the structure group of the covering $\ove{M}_i \to M$ 
%is isomorphic to $\zz$.
%The \ho~ $\xi_i: G \to \zz$
%%corresponding to the diagram \rrf{f:determ}
%sends the element $\t_i$ to
%a generator $\tau_i$ of $\zz$; we have 
%$\ove{f}_i(\tau_i x) = \ove{f}_i(x)-1$.
Fix a natural number $n$;
let 
$$
W_{n,i}
=
\bar f _i^{-1} \left([-n, 0]\right ) \sbs \ove M_i;
$$
this is a cobordism 
with two boundary components: $\bar f _i^{-1} (-n)$
and $\bar f _i^{-1} (0)$.
 Choose a $\tau_i$-ordered Morse function 
 $\Psi_i:
 W_{n,i}
 \to [ -1, m]$
 (where $m=\dim M$).
 It will be convenient to extend it to
 a function   
 $\ove M_i\to \rr$
 (denoted by the same symbol $\Psi_i$)
 in such  a way that 
$\Psi_i^{-1}\big([-1, m]\big) =  W_{n,i}$.
Let 
$$
\wh\Psi_i = \Psi_i\circ Q_i\circ Q; \ \ \ \ 
\Psi(x) = \max_i \wh\Psi_i(x); $$
\noindent 
then $\Psi$ is a continuous function on 
$\wh M$ and 
$$
\Psi^{-1}\left(]-\infty, m]\right)= Y, \ \ \ 
\Psi^{-1}\left(]-\infty, -1]\right)= Y_n.
$$
Put 
$$
Z_r = \{x~|~ \Psi(x) \leq r \}
\ \
{\rm (here } \ \ r= -1, \ldots, m).
$$
We will now  compute the homology of 
pairs 
$(Z_r,
Z_{r-1}
)$.
Recall that we use the following notation
$$
\lign =\big\{\l\in \zz G ~|~  \ \xi_i(\supp \l) 
\sbs 
]-\infty,-n]   \ {\rm for \ every }\  i   \big\},
$$
$$
\lgamin = \L_{\G, 0}^-; 
\ \ \ \
\L^{(n)}_{\G } = \lgamin\big/ \lign.
$$
Then $\lgamin$ is a ring with a unit
isomorphic to the  twisted
polynomial ring
$(\zz H)[\t_1, \ldots, \t_k]$
(where $H=\Ker\xi$)
and 
$
\L^{(n)}_{\G }
$
is the principal ideal of $\L_\G$
generated by the element $\t^n$
where $\t=\t_1\cdot\ldots \cdot \t_k$.
The ring 
$
\L^{(n)}_{\G }
$
is 
a free left $\zz H$-module generated
by monomials 
$\t^I$
where the multi-index $I\in \nn^k$ 
has at least one coordinate $\leq n-1$.

Since the functions $\Psi_i$ are $\tau_i$-ordered, 
the homology of every pair 
$(Z_r, Z_{r-1})$
is an $\L^{(n)}_{\G }$-module.
Let $\s$ be a 
zero of $p^*(\o)$ belonging to 
$Z_r \sm Z_{r-1}$.
 It is easy to see that the intersection $D_\s$ of the stable manifold of $\s$
with $Z_r\sm \Int Z_{r-1}$ is homeomorphic to an 
$r$-dimensional disc, and the pair $(D_\s, \pr D_\s)$ is homeomorphic
to $(D^r, S^{r-1})$. 
A zero of $p^*(\o)$ belonging to $Y\sm Y_n$
has index $r$ if and only if it is contained in $Z_r\sm Z_{r-1}$.
The chosen orientations of the stable manifolds of the flow generated
by $v$ determine  orientations of discs $D_\s$, compatible
with the action of the structure group of the covering $p$.
For every $\s$ the orientation of the disc $D_\s$ determines
 a homology class $\D_\s\in H_r(Z_r, Z_{r-1})$. Using the fact 
that the functions $\Psi_i$ are $t_i$-ordered
it is easy to see that 
\bq\lb{f:equiv}
\D_{\s g}=\D_\s \cdot g
{\rm \ \ 
for\  every  \ }
g\in \L^{(n)}_{\G }.
\end{equation}

\bele\lb{l:free_mod}
The $\L^{(n)}_{\G }$-module 
$F_r=
H_r(Z_r, Z_{r-1})$
is free with the basis
$\{\D_\s\}_{\s\in K_r}$. 
%$\qs$
\enle
\Prf
A standard Morse-theoretic 
argument implies that 
the homology of the 
pair $(Z_r, Z_{r-1})$
is non-zero only in degree $r$
and in this degree it is a free
abelian group generated by the 
elements 
$\D_{\s}$ where  
$\s$ ranges over the set of zeros of 
$p^*(\o)$  in 
$Z_r \sm Z_{r-1}$.
It is easy to see that every 
such zero $\s$ 
can be written uniquely 
in a form 
$a\cdot \t^I $ 
where 
 $a\in K_r$ 
and a multi-index $I=(i_1,\ldots, i_k)\in\nn^k$
has at least one coordinate $\leq n-1$. 
Applying the equality
\rrf{f:equiv} we complete the proof of the lemma. $\qs$

\bele\lb{l:same}
The map $\s\mapsto \D_\s$ 
determines a basis-preserving isomorphism of
the chain complexes $\NN_*^{(n)}$  and $F_*$.
\enle
\Prf
The assertion 
of the Lemma is equivalent to the statement that the
boundary operator of the 
homology exact sequence of the 
triple $(Z_r, Z_{r-1}, Z_{r-2})$
is given by the
formula \rrf{f:d_r} modulo $\t^n$.
%Fix some $i$ between $1$ and $n$.
We will deduce this  statement from the similar
one concerning 
the filtrations 
%induced by functions $\Psi_i$
discussed in Section \ref{s:circle}
(the case of circle-valued Morse functions).
Let
$$
X_i=\{x\in \wh M ~|~ \wh f_i(x)\leq 0\},\ \ \ \ \ \ \ 
X_{n,i}=\{x\in \wh M ~|~ \wh f_i(x)\leq -n\},
$$
so that 
$$Y=\bigcap_i X_i, \ \ \ \ \ Y_n = \bigcap_i X_{n,i}.$$

\noindent
Recall from Section \ref{s:circle}, \rrf{f:set-Y}, \rrf{f:filt-Y}
the filtration of $X_i$ induced by the function
$\wh\Psi_i$;
denote the terms of this filtration by $Z_{r,i}$ 
to emphasize here the dependence of these sets on 
$i\in \Lbrack 1, k \Rbrack$.
%(here $r\in \Lbrack 1, k \Rbrack$).
We have then 
$$
X_{n,i}= Z_{-1,i}\sbs Z_{0,i}\sbs \ldots Z_{m,i} = X_i.
$$

Denote by 
$\NN_{*,i}$
the  Novikov complex associated to the circle-valued Morse
function $f_i$, its gradient $v$ and the covering $p:\wh M \to M$.
Denote the corresponding truncated version by 
$\NN^{(n)}_{*,i}$.
We have an obvious commutative diagram

$$
\xymatrix{ 
\NN^{(n)}_{r}  \ar[r] \ar[d]^{A_i}  & H_r(Z_r, Z_{r-1}) \ar[d]^{B_i}\\
\NN^{(n)}_{r,i} \ar[r] & H_r(Z_{r,i}, Z_{r-1,i})
}
$$

Both horizontal arrows are isomorphisms, 
and the lower one commutes with the differentials
by Lemma \ref{l:isoo} of Section \ref{s:circle}.
Denote by $A$, resp. $B$ the direct sum of 
the maps $A_i$ and $B_i$. We  obtain the 
following diagram

$$
\xymatrix{ 
\NN^{(n)}_{r}  \ar[r] \ar[d]^{A}  & H_r(Z_r, Z_{r-1}) \ar[d]^{B}\\
\bigoplus_i \NN^{(n)}_{r,i} \ar[r] & \bigoplus_i H_r(Z_{r,i}, Z_{r-1,i})
}
$$
Both horizontal arrows are isomorphisms, the bottom arrow commutes 
with differentials, as well as the maps $A$ and  $B$.
It remains to observe that the map $A$ is injective by 
Corollary \ref{c:monomorphisms} and the proof of our Lemma is complete. $\qs$

The filtration of the chain complex $\SS_*(Y,Y_n)$ 
by subcomplexes $\SS_*(Z_r,Y_n)$ 
%(where $r\in\Lbrack 0,m \Rbrack$)
is therefore cellular. 
Similarly to the case of circle-valued Morse function
(see Lemma \ref{l:isoo} and the argument after it)
we conclude that there is a chain equivalence 
$$J_n:\NN_*^{(n)}\arrr \sim \SS_*(Y, Y_n)$$
as required by the first assertion of Proposition \ref{fl:kvadrat}.

Proceed to the proof of commutativity of the diagram
\rrf{f:kvadrat}.
Following \S 5 of \cite{P-o} (see p. 324, and Proposition A. 9) 
it is easy to prove that 
the chain equivalence $J_n$
%$\NN_*^{(n)} \to \SS_*(Y, Y_{n})$
does not depend up to homotopy 
on the choice of $\tau_i$-ordered 
Morse functions $\Psi_i$ on     $W_{n,i}$.
Let  $i$ be any integer in 
$ \Lbrack 1, k\Rbrack$.
It is shown  in \S  5 of \cite{P-o}
that for a given $n$ we can choose
$\tau_i$-ordered Morse functions
$\Psi_i:
 W_{n,i}
 \to [ -1, m]$
  and $\Phi_i : W_{n+1,i}
 \to [ -1, m]$
  in such a way that for every 
 $x\in  W_{n,i}$
 we have $\Phi_i(x)\geq \Psi_i(x)$,
 so that
 the inclusion of pairs 
 $(X_i, X_{n+1,i}) \arrinto (X_i, X_{n,i})$
 is compatible with the filtrations in these pairs induced 
 by $\Phi_i$, resp. $\Psi_i$.
Then the inclusion of pairs $(Y, Y_{n+1})\arrinto (Y, Y_n)$
 is also compatible with the filtrations in these pairs induced 
 by functions $\Phi = \max_i \Phi_i$ and $\Psi = \max_i \Psi_i$.
 The homotopic  commutativity of the square 
 \rrf{f:kvadrat} follows now from 
 Lemma 3.2 and Corollary 3.4 of \cite{P-o}.
The proof of Proposition \ref{fl:kvadrat}
is now over. $\qs$

% \NN_*^{(n)}

\subsection{Inverse limits of complexes and the end of the proof of Theorem \ref{t:main_th} }
\lb{su:inv_lims}

${}$

\noindent
We will now 
%use the Proposition 
%\ref{fl:kvadrat}
%to obtain 
construct
a chain equivalence between 
$\NN_*$ and $\steg$.
We have
\beeq\lb{f:limits}
\NN_*^- \approx \liminv \ \NN_*^{(n)}, 
\ \ \ \ \
\sximg \approx \ \liminv\ \SS_*(Y, Y_{n})
\end{equation}
(here $\FF=(\wh f_1, \dots, \wh f_k)$).
We wish to apply Proposition \ref{fl:kvadrat}
and build up the required chain equivalence
from the maps $J_n$.
However the diagrams 
\rrf{f:kvadrat}
commute only up to homotopy
and we need some more homological algebra to achieve the goal.
These algebraic tools will be developed in this section
(following largely  \cite{P-o}, \S\  3B).

\bele\lb{l:zil}(\cite{P-o}, Prop. 3.7. 1) )
Let
$$
A_*^1 \leftarrow \ldots \leftarrow A_*^s  \leftarrow \ldots,
\ \ \ \ \ 
B_*^1 \leftarrow \ldots \leftarrow B_*^s  \leftarrow \ldots
$$
be inverse sequences of chain complexes over a ring $R$.
Let $h_i:A_*^i \to 
B_*^i$ be
chain equivalences 
such that every square 
$$
\xymatrix{ 
A_*^i  \ar[d]^{h_i}   & A_*^{i+1} \ar[l] \ar[d]^{h_{i+1}}\\
B_*^i    & B_*^{i+1}\ar[l]
}
$$
\noindent
is homotopy commutative.
Let 
$$
\AAAA_*=\liminv \ A_*^i, \ \ \BBBB_* =\liminv\ B_*^i.
$$
Then there is a chain complex $\DDDD_*$ over $R$
and homology equivalences
$$ \AAAA_*
\arr
\DDDD_*
\arl 
\BBBB_*. \hspace{3cm} \qs$$

\enle
\bepr\lb{p:triangle}
Let 
$ \AAAA_*
,\
\DDDD_*
, \  
\BBBB_*$
be chain complexes over $R$
and 
$ \AAAA_*
\arrr \alpha
\DDDD_*
\arrl \beta 
\BBBB_*$
homology equivalences.
Assume that $\AAAA_*$ is a free chain complex.
Then there is a homology equivalence 
$\g: \AAAA_*\to \BBBB_*$ \sut~
$\b\circ\g\sim\a$.
\enpr

\Prf
First off we would like to replace the arrow $\b$
by an epimorphism.
\bede\lb{d:collaps}
Let $M$ be an $R$-module.
The chain complex
$$
\k_*=
\{
0
\arl
C_k
\arrl \pr
C_{k+1}
\arl
0
\}
$$
where 
$
C_k
=
C_{k+1}=M$
and $\pr=\id$
is called 
{\it elementary collapsible chain complex, concentrated in degrees 
$k, \ k+1$.}
A direct sum of elementary collapsible chain complexes
is called {\it collapsible chain complex}. $\qt$
\end{defi}
\bele\lb{l:making_epimorphic}
Let $\phi:\BBBB_*\to \FFFF_*$ be a chain map inducing 
an epimorphism in $H_0$. Then there is a 
collapsible chain complex $K_*$
%$\CCCC_*=\BBBB_*\oplus K_*$ 
%where $K_*$ is collapsible
and an epimorphic extension 
$\phi':\BBBB_*\oplus K_* \to \FFFF_* $ 
of the map $\phi$. 
\enle
\prf
Let 
$\k_*^k$ be the 
chain complex
$
\k_*=
\{
0
\leftarrow
\FFFF_k
\arrl \id 
\FFFF_{k}
\leftarrow
0
\}
$
concentrated in degrees 
$k-1, \ k$.
The maps
$\id: \FFFF_k\to \FFFF_k$
and 
$\pr: \FFFF_k\to \FFFF_{k-1}$
determine chain maps
$\m_k: \k_*^k \to \FFFF_*$.
Put $K_*=\bigoplus_{k\geq 1 } \k^k_*$.
We obtain a map 
$\mu: K_*\to \FFFF_*$.
Put 
%$\CCCC_*= \BBBB_*\oplus K_*$, and
$\phi'=(\phi, \mu)$.
The required properties are easy to check. $\qs$

\bele\lb{l:liftings}
Let 
$\a:\AAAA_*\to \FFFF_*$
and 
$\g:\EEEE_*\to \FFFF_*$
be homology equivalences of chain
complexes, \sut~ $\g$ is epimorphic and
$\AAAA_*$ is free. 
Then there is a chain map $\xi: \AAAA_*\to \EEEE_*$
\sut~ $\g\circ\xi=\a$.
\enle
\prf
The \ho s 
$\xi_r:\AAAA_r\to \FFFF_r$ 
commuting with boundary operators
are constructed by induction in $r$.
Assuming that we constructed $\xi_1, \ldots, \xi_r$
and
let $e$ be a free generator 
of $\AAAA_{r+1}$. 
Choose any element 
$y\in \EEEE_{r+1}$
\sut~
$\g(y)=\a(e)$.
Observe that $\z=\pr y -\xi_n(\pr e)$ 
is a cycle; moreover it is in the chain complex
$\ZZZZ_*=\Ker\g$, which is acyclic.
Therefore $\z=\pr \nu$
with $\nu\in \Ker\g$.
Put $\xi_{r+1}(e)=y-\nu$;
then 
$\xi_r(\pr e)=\pr \xi_{r+1}(e)$.
Similarly we extend the map $\xi_{r+1}$ to all the 
free generators
of $\AAAA_{r+1}$. $\qs$

\noindent
Proposition \ref{p:triangle}
is now easy to deduce. $\qs$ 

\pa
We can now complete the proof of Theorem \ref{t:main_th}.
Applying Lemma \ref{l:zil}
and Proposition \ref{p:triangle}
we obtain a homology equivalence
of $\Lgam$-complexes
$
\m: \nin\to\sximg$.
%Tensoring it by $\Lga$ 
Take the 
tensor product with $\Lga$, and
we obtain a chain map
$$
\wi \m = \mu \otimes \Id:
\nin\tens{\Lgam}\Lga  
\to
\sximg \tens{\Lgam}\Lga.
$$
\noindent
The proof of the next lemma is similar to proof of 
Proposition \ref{p:cones-tensor}
and will be omitted.
\bele\lb{l:mu}
The chain map $\wi\mu$ 
is a homology equivalence. $\qs$
\enle
\noindent
Applying Propositions \ref{p:cones-tensor}
and  Proposition \ref{p:compl-tens}
we obtain the following homology equivalences:
$$
\NN_*\approx\nin\tens{\Lgam}\Lga
\to 
\sxig
\leftarrow 
\steg.
$$
\noi
%Recall an isomorphism 
%$\sximg \tens{\Lgam}\Lga\approx \steg$ from Proposition \ref{p:cones-tensor}
%and a homology equivalence 
%$\steg \arrinto \sxig $
%from Proposition \ref{p:compl-tens}.
%that the natural inclusion 
%$\steg \arrinto \sxig $
%is a homology equivalence.
Applying again Proposition \ref{p:triangle}
we obtain a homology equivalence
\beeq\lb{f:final}
\NN_*\arrr \l \steg
\end{equation}
\noindent
Both chain complexes in \rrf{f:final}
are free, therefore $\l$ is a chain homotopy equivalence. 
The proof of Theorem \ref{t:main_th} is  over. $\qs$

%\newpage
%\section{White}
%\lb{s:White}

\newpage

\newpage
\section{Appendix 1. On the Pitcher inequalities for circle-valued Morse maps}
\label{s:appendix}

In the paper \cite{Pitcher}
Everett Pitcher obtained a lower bound for the number of 
critical points of a circle-valued 
Morse map. His remarkable 
work, dating back to 1939, is probably the first development 
in the circle-valued Morse theory.
In this Appendix we give an exposition of Pitcher's work
and relate it to the Novikov homology.
We show in particular that the numbers $Q_k$ introduced by Pitcher
equal Novikov Betti numbers, so that Pitcher's inequalities 
are equivalent to the torsion-free part of the Novikov inequalities.

\subsection{Pitcher inequalities}
\lb{su:pit_ineq}

We will use the terminology of E. Pitcher in order to stay as close as possible to his 
setup. Let $L$ be a closed \ma,
and $\t:L\to S^1=\rr/2\pi\zz$ a Morse map.
The image 
$
\t_*(H_1(L))\sbs H_1(S^1)=\zz
$
is a subgroup $\a\zz$ of $\zz$, where $\a$ is a positive integer. 
Let us assume for simplicity of exposition that $\a=1$.
Let $K\to L$ be the corresponding infinite cyclic covering,
and $F^*:K\to \rr$ be a Morse function making 
the following diagram commutative

$$
\xymatrix{ K    \ar[r]^{F^*}\ar[d] &  \rr \ar[d]
 \\
L\ar[r]^\t & S^1
}
$$

Denote by $T$ the generator of the structure group of
the covering, such that $F^*(Tx)=2\pi+F^*(x)$.
Choose a regular value $A$ of $F^*$
and let $B=A+2\pi$.
The set $\{x~|~F^*(x)\sbs [A, B[\ \}$
is a fundamental domain
for the action of the group $\zz$ on $K$.
To give the definition 
of the Pitcher's invariant 
$Q_k$, let us introduce some more terminology.
The closure of the fundamental domain above will be denoted by $W$;
it is a cobordism whose boundary is a disjoint union 
$(F^*)^{-1}(A) \sqcup (F^*)^{-1}(B)$ of two regular level surfaces of $F^*$.
Let $t=T^{-1}$, denote $(F^*)^{-1}(B)$ by $V$, then 
$\pr W= V\sqcup tV$. 
Let 
$V^-=(F^*)^{-1}(]-\infty, B])$
then $(F^*)^{-1}(]-\infty, A])= tV^-$.
E. Pitcher defines two numerical homological 
invariants of this configuration
(see the two paragraphs before Theorem I on the page 430 of \cite{Pitcher}).
\pa
1) {\it The count of the new $k$-cycles}.
\pa
In our notation this is the dimension
of the quotient
$H_k(V^-)/tH_k(V^-)$ (all homology groups are with rational coefficients).
Let us denote this number by $R_k$.
\pa
2) {\it The count of  newly bounding  $k$-cycles}.
\pa
This is the dimension
of $\Ker\Big(H_k(V^-) \arrr t H_k(V^-) \Big)$.
Let us denote this number by $S_k$.
The invariant $Q_k$ introduced by Pitcher is by definition
{\it the count of the new $k$-cycles
less
the count of  newly bounding  $k$-cycles},
that is,
\begin{equation}\lb{f:def_q}
 Q_k=R_k-S_k.
\end{equation}
We shall see a bit later, that $Q_k$ is a positive integer.
It follows immediately from the exact sequence of the pair
$(V^-, tV^-)$ that 
\begin{equation}\lb{f:betti}
b_k(V^-, tV^-)
=
R_k+S_{k-1}
\end{equation}
(where $b_k(V^-, tV^-)$ denotes the 
%$k$-th 
Betti number in degree $k$  of the pair 
$(V^-, tV^-)$).

Denote by $M_k$ the number of critical points of $\t$ of index $k$.
The classical Morse inequalities applied to
the function $F^*$ on the cobordism 
$W$
imply the inequalities
$$
M_k\geq b_k(V^-, tV^-) \geq R_k\geq Q_k.
$$
One deduces the
inequalities including the alternated sums of the above invariants:

\beth{ (Theorem I  \cite{Pitcher})}
for every $k$ we have
\bq\lb{f:pitch-ineq}
M_k-M_{k-1}+M_{k-2} +... \geq 
Q_k-Q_{k-1}+Q_{k-2} +... 
\end{equation}
\end{theo}

\prf 
Let us abbreviate $b_k(V^-, tV^-)$
to $\b_k$. Using the definition
\rrf{f:def_q}
of the Pitcher numbers $Q_k$ and the formula \rrf{f:betti}
it is easy to see that 
$$
\b_k-\b_{k-1}+\b_{k-2} +...
=
S_k
+
Q_k-Q_{k-1}+Q_{k-2} +... 
$$
The classical  Morse inequalities say

\bq\lb{f:m-ineq}
M_k-M_{k-1}+M_{k-2} +... \geq 
\b_k-\b_{k-1}+\b_{k-2} +...
\end{equation}
and the theorem follows. 
$\qs$ 

\subsection{Pitcher numbers and Novikov Betti numbers}
\lb{su:pit_nov}

Let $P=\qq[t]$.
The numbers $Q_k$ have a simple interpretation in terms
of the $P$-module structure 
of the homology $H_k(V^-)$.
The canonical decomposition of this module 
writes as follows:
\bq\lb{f:decomp}
H_k(V^-)
\approx
P^{a_k}
\oplus 
\big( \bigoplus_{i=1}^{c_k} P/t^{n_i}P\big)
\oplus 
\big( \bigoplus_{i=1}^{d_k} P/A_i P\big)
\end{equation}
where $n_i\in \nn, 0<n_i\leq n_{i+1}$
and $A_i\in \qq[t]$ 
are non-constant polynomials 
with non-zero free term, $A_i ~|~ A_{i+1}$.
We have then 
$$
R_k = a_k+c_k, \ \ \ S_k= c_k.
$$
Therefore $a_k=Q_k$.
Let $\L=\qq[t, t^{-1}]$.
Then
$H_*(K)\approx H_*(V^-)\tens{P}\L$,
and the rank of the $\L$-module $H_k(K)$
equals $a_k=Q_k$.
We deduce therefore the following
results:
\beth{(Theorem IV \cite{Pitcher} )}
\lb{t:hom_inv}
The numbers $Q_k$ are homotopy invariants of $L$ and the homotopy class of $\t$.
\enth

\beth 
The number $Q_k$ equals the Novikov Betti number 
$\wh b_k(L, [\t])$.
\enth

\newpage

\newpage
\section{Appendix 2. An example: Novikov complex whose incidence coefficient
has arbitrarily small convergence radius}
\label{s:appendixx}

\pa\pa

Let $q$ be any integer $\geq 3$.
In this Appendix we construct a circle-valued Morse function $f$
on a 3-manifold, and its transverse gradient
$u$ with the following properties:
\been\item The function $f$ has exactly two critical points
with indices 2 and  1. 
\item 
%There is a transverse $f$-gradient $u$
%\sut~ 
The unique  Novikov incidence coefficient  
is a power series of the form $\sum_k a_k t^k$
where $a_k\sim C\cdot q^k$,
and  $C\not=0$
(see the formula \rrf{f:inc-coef}).
\item 
This incidence coefficient  is stable \wrt~ 
$C^0$-small perturbations of the gradient.
\enen 

\noi The property 2) above implies that 
the convergence radius of the Novikov incidence coefficient
equals $1/q$, thus it converges to $0$ as $q\to\infty$.
The construction generalizes the example 
from the author's work \cite{Pepr2}, \S 3.
The proof uses the author's theory
of cellular gradients (\cite{P-stpet}, \cite{CVMT}).
and we begin by a brief outline of this theory.
The example itself is constructed in Section \ref{su:exam},
and the reader can start reading this section consulting 
the introductory Section \ref{s:cell_ratio}
when necessary.
\footnote{After this article was submitted to EJM,
the paper \cite{LM} of F. Laudenbach and C. Moraga appeared.
In this paper the authors announce a construction of a 
Morse-Novikov complex with infinite series coefficients.}

\subsection{Cellular gradients and rationality theorem}
\lb{s:cell_ratio}
${}$
\subsubsection{Cellular gradients of Morse functions on cobordisms }
\lb{su:cell_cobords}

${}$
\pa

Let $f:W\to[a,b]$ be a Morse function on a compact cobordism $W$;
put $\daw = f^{-1}(b), \ \dow = f^{-1}(a)$.
Pick an $f$-gradient $V$. For $a\in Crit(f)$ we denote by $D(a,v)$
the descending disc of $a$, that is, the stable manifold of $a$
\wrt~ flow induced by $v$. We denote by $D(v)$ the union of all descending discs
and by $D(\indl k; v)$ the union of all descending discs of critical 
points of indices $\leq k$. 
For $x\in \daw \sm D(-v)$
we denote by  $\stexp {(-v)}(x)$ the point where the $(-v)$-trajectory
$\g(x,t;-v)$ starting at $x$ intersects $\dow$.
The correspondence $x\mapsto \stexp {(-v)}(x)$ 
is then a diffeomorphism 
of $\daw \sm D(-v)$ onto $\dow \sm D(v)$. 
If $Crit(f)\not=\emptyset$ this map is not extensible to 
a continuous map of $\daw$ to $\dow$. However we have shown in 
\cite{P-stpet}, see also \cite{CVMT}, Part 3 
that for a $C^0$-generic gradient this  map can be endowed with a structure
that closely resembles a cellular map of a CW-complex. 

\bede\lb{d:MS-filtr}
\begin{itemize}\item
Let $\phi:N\to\rr$ be a self-indexing Morse function 
on a closed manifold $N$ (that is, $\phi~|~ Crit_k(\phi)=k$).
Put $N_i= \phi^{-1}\Big(]-\infty, i+1/2]\Big)$.
The filtration 
$$
\emptyset = N_{-1}\sbs N_0\sbs \ldots N_r =N$$
where $r=\dim N$ is called 
{\it the Morse-Smale filtration associated to $\phi$}
(or {\it MS-filtration} for brevity).
\item
For a given
MS-filtration $\{N_i\}$ of $N$
the filtration by submanifolds
$\wh N_j= \ove{N\sm N_{r-j-1}}$
is also an MS-filtration, called {\it the dual MS-filtration of the filtration $\{N_i\}$}
\end{itemize}
\end{defi}

\bere\lb{r:handles}
The term $N_s$ of an MS-filtration 
is the result of attaching to $N_{-1}$ of handles of indices $\leq s$;
it is a manifold with boundary 
homotopy equivalent to an $s$-dimensional CW-complex.
\end{rema}

\bede\lb{d:almost-transv}
An $f$-gradient 
$v$ is called {\it almost transverse}
if $D(p,v) \pitchfork D(q, -v)$ 
whenever
$\ind p \leq \ind q$.

The set of
 all $f$-gradients 
 is denoted by $G(f)$, the set of
 all almost transverse $f$-gradients 
is denoted by $G_A(f)$, the set 
of all  transverse $f$-gradients 
is denoted by $G_T(f)$.
\end{defi}

\bede\lb{d:cc}
\mfcob~ and
$v$ an almost transverse $f$-gradient.
We say that $v$ satisfies condition $(\gC)$
if there is a  Morse-Smale filtration
$\{\daw^k\}$ of $\daw$ and
a  Morse-Smale filtration
$\{\dow^k\}$ of $\dow$
such that for every $k$
\pa
\begin{align}\tag{$\gC$1}
\stv (\daw^k)
\sbs \hspace{.30cm} 
&
\Int \dow^k\hspace{.30cm}
\supset\hspace{.30cm}
D(\indl {k+1}, v)\cap \dow,
\lb{f:cc_1} \\
\tag{$\gC$2}
\st v \Big(~{\wh\dow~}^{k}~\Big)
\hspace{.30cm} \sbs \hspace{.30cm}
&
\Int {\wh\daw~}^{k}
\hspace{.30cm}
\supset\hspace{.30cm}
D(\indl {k+1}, -v)\cap \daw.\lb{f:cc_2}
\end{align}
\pa
The gradients satisfying condition $(\gC)$
will be also called {\it cellular gradients},
or $\gC$-gradients.
The set of all cellular gradients of $f$ will be denoted 
by $G_C(f)$.
\end{defi}

The following theorem  is one of the main results of \cite{P-stpet},
we cite it here using the terminology of \cite{CVMT}, Part 3.

\beth\lb{t:cell_dense}
The subset $G_C(f)\sbs G_A(f)$ is open and dense in $C^0$-topology.
\enth
Let $v$ be a cellular $f$-gradient for a Morse function $f$ on a cobordism $W$.  
Consider the compact topological space 
$\daw^k/\daw^{k-1}$ 
obtained by shrinking the subspace $\daw^{k-1}$ to a point
denoted $r_{k-1}$. 
The image of a point $y\in \daw^k$ in the space $\daw^k/\daw^{k-1}$
will be denoted by $\bar y$. 
Similar notation will be used 
for $\dow$, the  shrunk  subspace
$\dow^{k-1}$ will be denoted by $s_{k-1}$.
The next theorem describes  the cellular-like 
structure on the map $\stv$ 
(see \cite{CVMT}, p. 234).

\beth\lb{t:cell-map}
If $v$ is a cellular $f$-gradient then for every $k$ 
there is a continuous map 
$$\str {(-v)} : \daw^k/\daw^{k-1} \to \dow^k/\dow^{k-1}$$
\sut~ $\str {(-v)}(r_{k-1}) = s_{k-1}$ and 
\begin{align*}\lb{f:def_map}
&   \str {(-v)}(\ove x)\ =\ s_{k-1} &
\ 
\mxx{ if } &\ \ 
x\in D(-v);             \\
&\str {(-v)} (\ove x) =\ove{\stexp {(-v)}(x)}
& 
\mxx{ if } & x\notin  D(-v).\ \ \ 
%\qs
\end{align*}
\enth

\bede\lb{d:homol_grad_desc}
The map induced by $\str {(-v)}$
in homology is denoted by
\bqq
H_k(-v):
H_k(\daw^k, \daw^{k-1}  )
\to
H_k(\dow^k, \dow^{k-1})
\end{equation*}
\noi
and  called {\it homological gradient descent}.
\end{defi}
This \ho~ is stable \wrt~ $C^0$-small perturbations of the gradient, as shown
in \cite{CVMT}, p.282:

\bepr\lb{p:pert_homol}
Let $v$ be a cellular gradient of a Morse function
$f:W\to [a,b]$.
There is $\d>0$ such that for every
$f$-gradient $w$ with $||w-v||<\d$
the homomorphisms
$$H_k(-v),~ H_k(-w)  :
H_k(\daw^k, \daw^{k-1})
\to
H_k(\dow^k, \dow^{k-1})
$$
are equal. 
%$\qs$
\enpr

\subsubsection{Cellular gradients for circle-valued Morse functions
and their Novikov complexes}
\lb{su:cell_circle}

${}$
\pa

Let $f:M\to S^1=\rr/\zz$ be a Morse function,
we will assume that its class $[f]$ in $H^1(M,\zz)$ 
is indivisible. 
%Let  $\l$ be a  regular value of $f$.
Let $\ove M \to M$ be the corresponding infinite cyclic covering;
lift the function $f$ to a real-valued Morse function $F:\ove M\to\rr$.
Let $\l$ be a regular value of $F$,
put $V=F^{-1}(\l)$. We have a \co~ $W=F^{-1}([\l-1,\l])$
and  a Morse function $F|W:W\to [\l-1, \l]$.
Let $t$ be a generator of the structure group $\approx \zz$ of the covering,
\sut~ $F(tx)=F(x)-1$. The 
map $t^{-1}$ determines a diffeomorphism $\dow\to\daw$
which will  be denoted by $I$. 
%Cutting $M$ along $f^{-1}(\l)$ 
%we obtain a \co~ $W$, a Morse function $F:W\to [0,1]$
%and a diffeomorphism $I:\daw\to\dow$.
An $f$-gradient $v$ induces an $F$-gradient, denoted
by the same symbol $v$.

\bede\lb{d:cell_circle}
\begin{itemize}
 \item 
 An $f$-gradient $v$ is called {\it cellular \wrt~ $\l$} if 
 the induced $F$-gradient on $W$ is cellular \wrt~ 
 some MS-filtration $\{N_i\}$ on $\dow$ and 
 the MS-filtration $\{I(N_i)\}$ on $\daw$.
 \item An $f$-gradient $v$ is called {\it cellular} if 
 it is cellular \wrt~ $\l$ for some regular value $\l$ of $f$.
 \item The 
 set of all cellular gradients of $f$ is denoted $G_C(f)$.
 \end{itemize}
 \end{defi}
 
 The following theorem  is one of the main results of \cite{P-stpet},
 concerning circle-valued Morse functions;
we cite it here in the terminology of \cite{CVMT}, Ch. 12.
 
\beth\lb{t:circle-cell-dense}
\been\item The subset $G_C(f)\sbs G(f)$ 
is open and dense in $G(f)$ \wrt~ 
$C^0$-topology.
\item 
The subset $G_C(f)\cap G_T(f)\sbs G_T(f)$ 
is open and dense in $G_T(f)$ \wrt~ 
$C^0$-topology.
\enen
\enth

Let $v$ be a cellular $f$-gradient.
For every $k$ we have an endomorphism
 \bqq
\HH_k(-v)= I_*\circ H_k(-v):
H_k(\daw^k, \daw^{k-1}  )
\to
H_k(\daw^k, \daw^{k-1}  ).
\end{equation*}

\noindent
The Proposition \ref{p:pert_homol}
implies the following Corollary.

\beco\lb{c:cell_perturb}
Let $v$ be a cellular $f$-gradient. There is $\d>0$
\sut~ for every $f$-gradient $w$ 
with $||v-w||<\e$ 
and every $r$ we have 
$\HH_k(-v)=\HH_k(-w)$.
\enco 
\noi
It turns out that 
the Novikov complex $N_*(f,v)$ 
can be computed in terms of this \ho. 
Let $p\in Crit_{k+1}(f),
q\in Crit_k(f)$.
Choose the lifts $\bar p, \bar q$ of the points $p, q$ to $\ove M$ 
in such a way that $\bar p, \ \bar q\in t^{-1}W$.
Since $v$ is cellular, 
the $k$-dimensional \sma~ 
$T=D(\bar p,v)\cap \daw$ is in $\daw^k$,
and the set 
$\daw^k\sm \ \Int \daw^{k-1}$ is compact. 
Choose some orientations of descending discs $D(p,v), \ D(q,v)$.
We have then the fundamental class $[T]\in H_k(\daw^k, \daw^{k-1})$.
Assume for simplicity of exposition that $M$ is oriented.
Similarly, the $(r-k)$-dimensional \sma~ 
$S=D(t\bar q, -v)\cap \daw$
determines a class
$$[S]\in H_{r-k}\Big(\wh\daw^{r-k}, \wh\daw^{r-k-1}\Big),$$
where $r=\dim V=\dim M -1$.
The intersection index $\langle [T], [S] \rangle\in \zz$
is defined.
%Denote by $A$ the endomorphism $\HH_k(-v)$.
The next theorem (see \cite{CVMT}, p. 379)
expresses the 
% %Here $N(p,q;v)$ denotes the
Novikov incidence coefficient 
$N(p,q;v)\in \zz((t))$
in terms of the gradient descent \ho.

\beth\lb{t:comput}
We have
\bqq
N(p,q;v)
=
n_0(p,q;v) + \sum_{m\geq 0} \Big\langle (\HH_k(-v))^m\big([T]\big), [S]\Big\rangle
t^{m+1}.
%\in \zz[[t]].
\end{equation*}
\enth
\noi
(Here $n_0(p,q;v)\in\zz$ is the incidence coefficient of the critical 
points $p, q$ in the \co~ $W$;
the brackets $\langle \cdot\ , \  \cdot\rangle$ 
denote the intersection index.) 

The next Corollary is obtained by a standard argument from linear algebra.
\beco\lb{c:ratio}
For any cellular gradient $v$ the 
 Novikov incidence coefficient $N(p,q;v)$ is a rational
function of the form
$\frac {P(t)}{Q(t)}$
where $P,Q\in\zz[t]$ and $Q(0)=1$.
\enco

%$\st {w_1}$

%$\stind a12$

\subsection{An example}
\lb{su:exam}

${}$
\pa
Let $\ttt^2$ be the 2-dimensional torus, $\alpha$ be its parallel,
$\beta$ its meridian.
Consider two disjoint closed discs $D_1, D_2$ in $\ttt^2$
which do no intersect $\alpha\cup\b$.
Removing their interiors from $\ttt^2$ we obtain 
a surface $S$, whose boundary is the disjoint union of
two circles $\pr_1 S$ and $\pr_2 S$
(see the upper image on the Fig. 1).

Attach a copy $S(1,1)$ of 
$S$ to another copy $S(1,2)$ of $S$,
identifying 
$\pr_2 S(1,1)$
with 
$\pr_1 S(1,2)$.
We obtain a surface of genus 2 
with two components of boundary:
$\pr_1 S(1,1)$
and
$\pr_2 S(1,2)$.
Attaching the copies $D_1(1), D_2(1)$ of the discs
$D_1, D_2$ to these components
gives a closed surface $N$ of genus 2. 
One more copy of this surface will be denoted by $K$
(see the bottom of the Fig. 1).
Similarly we glue together 
three copies 
$S(1/2,1),  S(1/2,0), S(1/2,2)$
of $S$ and attach 
to it two discs $D_1(1/2), D_2(1/2)$ 
to obtain a closed surface $L$ of genus 3 
(depicted in the middle of the figure).
Associate to every  point in 
$ D_1(1/2)\cup
S(1/2,1)
\cup S(1/2,0)
$ 
its copy in 
$
D_1(1)\cup
S(1,1)
\cup S(1,2)
$;
this determines a diffeomorphism which
will be denoted by 
$I(1/2, 1)$.
Similarly, we construct a diffeomorphism
$I(1/2, 0)$
of the surface
$
S(1/2,0)\cup S(1/2,2)\cup 
D_1(1/2)$
onto
$S(0,1)\cup S(0,2)
\cup D_2(0).
$

A surgery along the circle
$\b(1/2, 2)$ 
yields 
a surface naturally diffeomorphic to $N$.
Attaching the corresponding handle of index 2
to $L\times [0, \e]$ gives a \co~ $W_1$
endowed with a Morse function $F_1:W_1\to [1/2, 1]$.
This Morse function has one critical point $x_2$ of index 2.
Pick a gradient $w_1$ for this function in such a way
that 
\been\item[A1) ]
The ascending disc $D(x_2, -w_1)$ intersects
the level surface $N=F^{-1}_1(1)$ by two points in the 
interior of $D_2(1)$.
\item[A2) ] The diffeomorphism
$\st {w_1}$ 
sends $D_2(1/2)$ to the interior of $D_2(1)$.
\item[A3) ] 
The restriction of 
$\st {w_1}$ 
to
$D_1(1/2)\cup S(1/2,1)\cup S(1/2,0)$
equals 
$I(1/2, 1)$
everywhere except a small tubular \nei~ $T_1$ of the circle 
$\pr D_1(1/2)$.
Further,
$\st {w_1} (T_1) = I(1/2, 1) (T_1)$
and the $\st {w_1}$-image of 
 $D_1(1/2)$
contains 
$D_1(1)$
in its interior.
\enen

Similarly, we do a surgery along the circle 
$\b(1/2,1)$
and obtain a surface naturally diffeomorphic to $K$.
Attach the corresponding handle to $L\times [-\e, 0]$,
get a \co~ $W_0$ endowed with a Morse function 
$F_0: W_0\to [0,1/2]$ having one critical point $x_1$ of index 1.
Pick an $F_0$-gradient $w_0$ \sut~

\been\item[B1) ]
The descending disc $D(x_1, w_0)$ intersects
the level surface $K=F^{-1}_0(0)$ by two points in the 
interior of $D_1(0)$.
\item[B2) ] The diffeomorphism
$\stexp {(-w_0)}$ 
sends $D_1(1/2)$ to the interior of $D_1(0)$.
\item[B3) ] 
The restriction of 
$\stexp {(-w_0)}$ 
to
$ S(1/2,0)\cup S(1/2,2)\cup D_1(1/2)$
equals 
$I(1/2, 0)$
everywhere except a small tubular \nei~ $T_2$ of the circle 
$\pr D_2(1/2)$.
Further,
$\stexp {(-w_0)}(T_2) = I(1/2, 0) (T_2)$
and the $\stexp {(-w_0)}$-image of 
 $D_2(1/2)$
contains 
$D_2(0)$
in its interior.
\enen
We have $\pr W_1 \approx N\sqcup L, \ \pr W_0 \approx L\sqcup K$;
attaching $W_1 $ to $W_0$ 
along the $L$-component  of their boundaries 
%which is 
%diffeomorphic to $ L$
we obtain a \co~ $W$ with boundary $\pr W\approx N\sqcup K$,
endowed with a Morse function $F:W\to [0,1]$,
\sut~ 
$Crit(F)=\{x_1, x_2\}, \ \ind x_j = j$.
The gradients $w_0$ and $w_1$ can be glued together 
(modifying them appropriately nearby $L$ if necessary)
so that the resulting gradient $w$ is  cellular $F$-gradient
(see Definition \ref{d:cc}). To show this we introduce   
Morse-Smale filtrations on $\dow$ and $\daw$. 
Let 
$$N_0=D_1(1),
\ \ N_1 = 
  D_1(1) \cup S(1,1)\cup S(1,2),
  \ \ 
  N_2= N.
  $$
  The filtration $N_0\sbs N_1\sbs N_2$
  is then a MS-filtration of $N$.
  The image of this filtration 
  \wrt~ the natural diffeomorphism 
  $J: N \to K$
  is a  MS-filtration $K_0\sbs K_1 \ \sbs K_2$ on $K$.
    The properties A1) -- A3) and B1) -- B3) 
  imply the following:
  
  \begin{align*}
   & (D1) &\ D(x_i, w) \cap \dow \sbs \Int K_{i-1} \ \ &{\rm for } \ \ i=1,2  \\
   & (D2) &\stexp {(-w)} (N_i) \sbs \Int K_i \ \ &{\rm for } \ \ i=0,1,2.
  \end{align*} 
  \noindent
  And we have also the dual properties
  \noindent
  \begin{align*}
   & (U1) &\ D(x_i, -w) \cap \daw \sbs \Int \wh N_{2-i} \ \ &{\rm for } \ \ i=1,2  \\
   & (U2) &\st {w} (\wh K_i) \sbs \Int \wh N_i \ \ &{\rm for } \ \ i=0,1,2.
  \end{align*} 
  \noi
  (recall that $\wh N_i$ and $\wh K_i$ denote the MS-filtrations
  dual to the filtrations 
$N_i$, resp.  $ K_i$).
The conjunction of the properties $(D1), (D2)$ is just a reformulation of the 
condition \rrf{f:cc_1}; similarly the conjunction of the properties $(U1), (U2)$
is  equivalent to \rrf{f:cc_2}. The $F$-gradient $w$ is therefore cellular.

The 3-manifold 
$M$ and a circle-valued function on it will be obtained by gluing 
$N$ to $K$ via a diffeomorphism that we will now describe.
Put
$$
a_1=[\alpha(1,1)], \
a_2=[\alpha(1,2)], \
b_1=[\b(1,1)], \ 
b_2=[\b(1,2)]. \ 
$$
The family $\BB=(a_1, b_1, a_2, b_2)$
is then a basis in $H_1(N)$.
The same embedded circles determine the homology classes in 
$H_1(N_1/N_0)$, they will be denoted by the same letters by a certain
abuse of notation.
Similarly we obtain a base $\CC=(a'_1, b'_1, a'_2, b'_2)$
in $H_1(K)$ and $H_1(K_1/K_0)$. Denote by $J:N\arrr \approx K$ 
the natural diffeomorphism, then we have $J_*(\BB)=\CC$.
Consider the automorphism of $H_1(K)\approx \zz^4$
given in the base $\CC$ by the following  matrix

$$
\sujap=\left( \begin{matrix} 0 & 2 & 1 & 0\\
0 & 0 & 0 & 1\\
0 & 1 & 0 & 0\\
-1 & q & 0 & -2\\
\end{matrix}
\right)
$$
where $q$ is any  integer $\geq 3$. 
It is easy to check that $\sujap$ preserves the intersection form on $H_1(N)$
therefore there is a diffeomorphism $\Phi:K\to K$
inducing $\sujap$ in $H_1$. We can assume that
$\Phi(x)=x$ for $x\in D_1(0)\cup D_2(0)$. 
Identifying each point $y\in K$ with $J\Phi(y)\in N$
%Gluing $K$ to $N$ via the diffeomorphism $J^{-1}\circ \Phi$
we obtain a 3-manifold $M$; the 
Morse function $F$ induces a map $f:M\to S^1$.
The $F$-gradient $w$ induces an $f$-gradient $v$.
Observe that $v$ is a cellular $f$-gradient \wrt~
regular level surface $N$ and its MS-filtration $\{N_i\}$.
The matrix of the  endomorphism  
$\HH_1(-v)~:~H_1(N_1/N_0)\to H_1(N_1/N_0)$
is easy to compute; it equals
$$
\sajap=
\left(
\begin{matrix}
0 & 0 & 0 & 2\\
0 & 0 & 0 & 0\\
0 & 0 & 0 & 1\\
0 & 0 & -1 & q
\end{matrix}
\right).
$$
Pick a transverse $f$-gradient $u$ sufficiently  close to $v$ in $\smo$ 
topology so
that $u$ is still a cellular $f$-gradient \wrt~ the level surface $N$
and its MS-filtration $\{N_i\}$, and $\HH_1(-u)=\HH_1(-v)$.
The Novikov incidence coefficient
$N( x_2,  x_1; u)$
is now  easy to compute.
Let $T$ be the $J^{-1}$-image in  $N$ of 
$D(x_2,u)\cap K$, and $\t=[T]\in H_1(N_1/N_0)$.
Then $\t=b_1-2b_2$.
Let $S=D(x_1,-u)\cap N$ 
%be the image in  $M$ of 
%$D(x_1,-u)\cap N$, 
then  $[S]=b_1\in H_1(N_1/N_0)$.
Applying Theorem \ref{t:comput} 
we obtain
$$n_{k+1}( x_2,  x_1; u)
=
\big\langle \sajap^k(\t), b_1\big\rangle.
$$
We have 
$\sajap^k(\t)=(-2)\cdot  \sajap^k(b_2)$. 
Therefore 
$n_{k+1}( x_2,  x_1; u)$
equals the first coordinate of the vector 
$(-2)\cdot \sajap^k(b_2)$ \wrt~ basis $\BB$.
Computing this first coordinate 
is a routine exercise in linear algebra 
which will be left to the reader.
We give just the result:
\begin{gather}\lb{f:inc-coef}
  n_{k+1}( x_2,  x_1; u)
 =
 \frac {-4}{\sqrt{q^2-4}} \Big(A^k - B^k\Big), \\
 {\rm where } ~~ 
 A=(q+\sqrt{q^2-4})/2; \ \ \ \  B = 
  (q-\sqrt{q^2-4})/2.
\end{gather}
\noindent
The properties 
of $f$ and $v$ stated 
in the beginning of this Appendix 
are now obvious. 
${}$
\newpage 
$$
\hspace{-4cm}
\includegraphics[height=20cm]{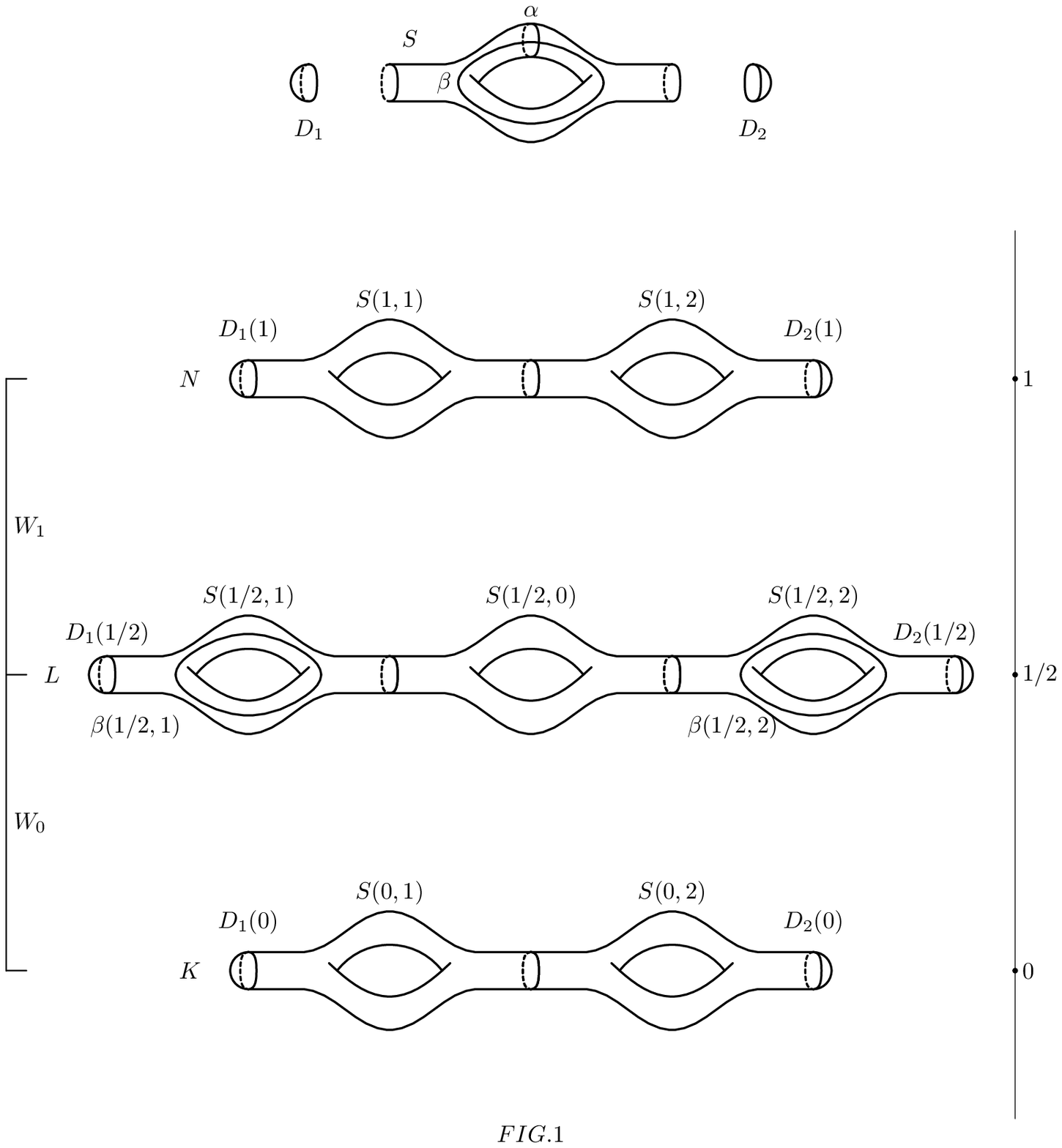}
\vspace{10cm}
$$

\newpage

\newpage

\bibliography{mybib}
\bibliographystyle{plain}

\end{document}